\documentclass{gtart}

\def\ifplaintex{\expandafter\ifx\csname documentclass\endcsname\relax}

\def\gtp{{\mathsurround=0pt\it $\cal G\mskip-2mu$eometry \&\ 
$\cal T\!\!$opology $\cal P\!$ublications}}  

\def\recd{{\small Received:\qua\receiveddate\ifx\reviseddate\relax
\else\qquad Revised:\qua\reviseddate\fi\par}} 

\def\lognumber#1{\def\thelognumber{#1}}
\def\volumenumber#1{\def\thevolumenumber{#1}}
\def\volumeyear#1{\def\thevolumeyear{#1}}
\def\papernumber#1{\def\thepapernumber{#1}}
\def\pagenumbers#1#2{\def\startpage{#1}\def\finishpage{#2}}
\def\published#1{\def\publishdate{#1}}

\def\received#1{\def\receiveddate{#1}}
\def\revised#1{\def\reviseddate{#1}}
\def\accepted#1{\def\accepteddate{#1}}
\def\asciititle#1{\def\theasciititle{#1}}

\long\def\asciiabstract#1{\long\def\theasciiabstract{#1}}
\def\asciikeywords#1{\def\theasciikeywords{#1}}

\let\\\par\let\thelognumber\relax\let\thevolumenumber\relax
\let\thepapernumber\relax\let\thevolumeyear\relax\let\startpage\relax
\let\finishpage\relax\let\publishdate\relax\let\receiveddate\relax
\let\reviseddate\relax\let\accepteddate\relax\let\theasciititle\relax
\let\theasciiauthors\relax
\let\theasciiabstract\relax\let\theasciikeywords\relax

\let\theasciiemail\relax

\ifplaintex
\font\logobig=cmssbx10 scaled 3836
\font\logomed=cmssbx10 scaled 2557
\else
\font\logobig=cmssbx10 scaled 4200
\font\logomed=cmssbx10 scaled 2800
\fi

\long\def\makeagttitle{   
\count0=\startpage
\agt\hfill      
\hbox to 45truept{\vbox to 0pt{\vglue -13truept{\logomed A\kern -.37em{\logobig 
T}\kern -.38em G}\vss}\hss}
\break
{\small Volume \thevolumenumber\ (\thevolumeyear)
\startpage--\finishpage\nl
Published: \publishdate}

\vglue .25truein

{\parskip=0pt\leftskip 0pt plus
1fil\def\\{\par\smallskip}{\Large\bf\thetitle}\par\medskip} \vglue
0.05truein

{\parskip=0pt\leftskip 0pt plus 1fil\def\\{\par}{\sc\theauthors}
\par\medskip}%
 
\vglue 0.03truein 

{\small\leftskip 25truept\rightskip 25truept{\bf Abstract}\stdspace\theabstract

{\bf AMS Classification}\stdspace\theprimaryclass
\ifx\thesecondaryclass\relax\else; \thesecondaryclass\fi\par
{\bf Keywords}\stdspace \thekeywords\par}\vglue 7truept

}   

\ifplaintex
\font\phead=cmsl9 scaled 950
\font\pnum=cmbx10 scaled 913
\font\pfoot=cmsl9 scaled 950
\headline{\vbox to 0pt{\vskip -4.5mm\line{\small\phead\ifnum
\count0=\startpage ISSN 1472-2739 (on-line) 1472-2747 (printed)
\hfill {\pnum\folio}\else\ifodd\count0\def\\{ }%
\ifx\theshorttitle\relax\thetitle\else\theshorttitle\fi\hfill{\pnum\folio}
\else\def\\{ and }{\pnum\folio}\hfill\ifx\theshortauthors\relax\theauthors
\else\theshortauthors\fi\fi\fi}\vss}}
\footline{\vbox to 0pt{\vglue 0mm\line{\small\pfoot\ifnum\count0=\startpage
\copyright\ \gtp\hfill\else
\agt, Volume \thevolumenumber\ (\thevolumeyear)\hfill\fi}\vss}}
\else
\font=cmsl9 scaled 1050
\font\lnum=cmbx10 
\font=cmsl9 scaled 1050
\makeatletter
\def\@oddhead{{\small\lhead\ifnum\count0=\startpage ISSN 1472-2739 
(on-line) 1472-2747 (printed)\hfill {\lnum\number\count0}\else\ifodd\count0
\def\\{ }\ifx\theshorttitle\relax \thetitle \else\theshorttitle\fi\hfill
{\lnum\number\count0}\else\def\\{ and }{\lnum\number\count0}
\hfill\ifx\theshortauthors\relax 
\theauthors\else\theshortauthors\fi\fi\fi}}\def\@evenhead{\@oddhead}
\def\@oddfoot{\small\lfoot\ifnum\count0=\startpage\copyright\ \gtp\hfill\else
\agt, Volume \thevolumenumber\ (\thevolumeyear)\hfill\fi}
\def\@evenfoot{\@oddfoot}
\makeatother
\fi
\let\maketitlepage\makeagttitle
\let\makeshorttitle\maketitlepage
\let\maketitle\maketitlepage

\newwrite\gtoutfile
\long\gdef\makeheadfile{  
{\def\\{, }\def\s{ }
\immediate\openout\gtoutfile head.xxx
\immediate\write\gtoutfile{To: math@arxiv.org}
\immediate\write\gtoutfile{Subject: put OR rep NNNNN:ppppp}
\immediate\write\gtoutfile{--text follows this line--}
\immediate\write\gtoutfile{Proxy-for: \ifx\theasciiauthors\relax
\theauthors\else\theasciiauthors\fi\s<\ifx\theasciiemail\relax\theemail\else\theasciiemail\fi>}
\immediate\write\gtoutfile{\noexpand\\}
\immediate\write\gtoutfile{Authors: \ifx\theasciiauthors\relax
\theauthors\else\theasciiauthors\fi}
{\def\\{ }\immediate\write\gtoutfile{Title: \ifx\theasciititle\relax
\thetitle\else\theasciititle\fi}}
\immediate\write\gtoutfile{Subj-class: GT or SG, GR etc}
\immediate\write\gtoutfile{MSC-class: \theprimaryclass\ifx\thesecondaryclass\relax\else, \thesecondaryclass\fi}
\immediate\write\gtoutfile{Journal-ref: Algebr. Geom. Topol. \thevolumenumber\s
(\thevolumeyear) \startpage-\finishpage}
\immediate\write\gtoutfile{Comments: Published by Algebraic and
Geometric Topology at}
\immediate\write\gtoutfile{\s\s\s  http://www.maths.warwick.ac.uk/agt/AGTVol\thevolumenumber/agt-\thevolumenumber-\thepapernumber.abs.html}
\immediate\write\gtoutfile{\noexpand\\}
\immediate\write\gtoutfile{}
\ifx\theasciiabstract\relax
\immediate\write\gtoutfile{\theabstract}\else
\immediate\write\gtoutfile{\theasciiabstract}\fi
\immediate\write\gtoutfile{}
\immediate\write\gtoutfile{\noexpand\\}
\immediate\write\gtoutfile{}
\immediate\closeout\gtoutfile}}  

\def\maketitlepage{\makeagttitle\makeheadfile}
\let\makeshorttitle\maketitlepage
\let\maketitle\maketitlepage

\def\ifplaintex{\expandafter\ifx\csname documentclass\endcsname\relax}

\def\gtp{{\mathsurround=0pt\it $\cal G\mskip-2mu$eometry \&\ 
$\cal T\!\!$opology $\cal P\!$ublications}}  

\def\recd{{\small Received:\qua\receiveddate\ifx\reviseddate\relax
\else\qquad Revised:\qua\reviseddate\fi\par}} 

\def\lognumber#1{\def\thelognumber{#1}}
\def\volumenumber#1{\def\thevolumenumber{#1}}
\def\volumeyear#1{\def\thevolumeyear{#1}}
\def\papernumber#1{\def\thepapernumber{#1}}
\def\pagenumbers#1#2{\def\startpage{#1}\def\finishpage{#2}}
\def\published#1{\def\publishdate{#1}}

\def\received#1{\def\receiveddate{#1}}
\def\revised#1{\def\reviseddate{#1}}
\def\accepted#1{\def\accepteddate{#1}}
\def\asciititle#1{\def\theasciititle{#1}}

\long\def\asciiabstract#1{\long\def\theasciiabstract{#1}}
\def\asciikeywords#1{\def\theasciikeywords{#1}}

\let\\\par\let\thelognumber\relax\let\thevolumenumber\relax
\let\thepapernumber\relax\let\thevolumeyear\relax\let\startpage\relax
\let\finishpage\relax\let\publishdate\relax\let\receiveddate\relax
\let\reviseddate\relax\let\accepteddate\relax\let\theasciititle\relax
\let\theasciiauthors\relax
\let\theasciiabstract\relax\let\theasciikeywords\relax

\let\theasciiemail\relax

\ifplaintex
\font\logobig=cmssbx10 scaled 3836
\font\logomed=cmssbx10 scaled 2557
\else
\font\logobig=cmssbx10 scaled 4200
\font\logomed=cmssbx10 scaled 2800
\fi

\long\def\makeagttitle{   
\count0=\startpage
\agt\hfill      
\hbox to 45truept{\vbox to 0pt{\vglue -13truept{\logomed A\kern -.37em{\logobig 
T}\kern -.38em G}\vss}\hss}
\break
{\small Volume \thevolumenumber\ (\thevolumeyear)
\startpage--\finishpage\nl
Published: \publishdate}

\vglue .25truein

{\parskip=0pt\leftskip 0pt plus
1fil\def\\{\par\smallskip}{\Large\bf\thetitle}\par\medskip} \vglue
0.05truein

{\parskip=0pt\leftskip 0pt plus 1fil\def\\{\par}{\sc\theauthors}
\par\medskip}%
 
\vglue 0.03truein 

{\small\leftskip 25truept\rightskip 25truept{\bf Abstract}\stdspace\theabstract

{\bf AMS Classification}\stdspace\theprimaryclass
\ifx\thesecondaryclass\relax\else; \thesecondaryclass\fi\par
{\bf Keywords}\stdspace \thekeywords\par}\vglue 7truept

}   

\ifplaintex
\font\phead=cmsl9 scaled 950
\font\pnum=cmbx10 scaled 913
\font\pfoot=cmsl9 scaled 950
\headline{\vbox to 0pt{\vskip -4.5mm\line{\small\phead\ifnum
\count0=\startpage ISSN 1472-2739 (on-line) 1472-2747 (printed)
\hfill {\pnum\folio}\else\ifodd\count0\def\\{ }%
\ifx\theshorttitle\relax\thetitle\else\theshorttitle\fi\hfill{\pnum\folio}
\else\def\\{ and }{\pnum\folio}\hfill\ifx\theshortauthors\relax\theauthors
\else\theshortauthors\fi\fi\fi}\vss}}
\footline{\vbox to 0pt{\vglue 0mm\line{\small\pfoot\ifnum\count0=\startpage
\copyright\ \gtp\hfill\else
\agt, Volume \thevolumenumber\ (\thevolumeyear)\hfill\fi}\vss}}
\else
\font=cmsl9 scaled 1050
\font\lnum=cmbx10 
\font=cmsl9 scaled 1050
\makeatletter
\def\@oddhead{{\small\lhead\ifnum\count0=\startpage ISSN 1472-2739 
(on-line) 1472-2747 (printed)\hfill {\lnum\number\count0}\else\ifodd\count0
\def\\{ }\ifx\theshorttitle\relax \thetitle \else\theshorttitle\fi\hfill
{\lnum\number\count0}\else\def\\{ and }{\lnum\number\count0}
\hfill\ifx\theshortauthors\relax 
\theauthors\else\theshortauthors\fi\fi\fi}}\def\@evenhead{\@oddhead}
\def\@oddfoot{\small\lfoot\ifnum\count0=\startpage\copyright\ \gtp\hfill\else
\agt, Volume \thevolumenumber\ (\thevolumeyear)\hfill\fi}
\def\@evenfoot{\@oddfoot}
\makeatother
\fi
\let\maketitlepage\makeagttitle
\let\makeshorttitle\maketitlepage
\let\maketitle\maketitlepage

\newwrite\gtoutfile
\long\gdef\makeheadfile{  
{\def\\{, }\def\s{ }
\immediate\openout\gtoutfile head.xxx
\immediate\write\gtoutfile{To: math@arxiv.org}
\immediate\write\gtoutfile{Subject: put OR rep NNNNN:ppppp}
\immediate\write\gtoutfile{--text follows this line--}
\immediate\write\gtoutfile{Proxy-for: \ifx\theasciiauthors\relax
\theauthors\else\theasciiauthors\fi\s<\ifx\theasciiemail\relax\theemail\else\theasciiemail\fi>}
\immediate\write\gtoutfile{\noexpand\\}
\immediate\write\gtoutfile{Authors: \ifx\theasciiauthors\relax
\theauthors\else\theasciiauthors\fi}
{\def\\{ }\immediate\write\gtoutfile{Title: \ifx\theasciititle\relax
\thetitle\else\theasciititle\fi}}
\immediate\write\gtoutfile{Subj-class: GT or SG, GR etc}
\immediate\write\gtoutfile{MSC-class: \theprimaryclass\ifx\thesecondaryclass\relax\else, \thesecondaryclass\fi}
\immediate\write\gtoutfile{Journal-ref: Algebr. Geom. Topol. \thevolumenumber\s
(\thevolumeyear) \startpage-\finishpage}
\immediate\write\gtoutfile{Comments: Published by Algebraic and
Geometric Topology at}
\immediate\write\gtoutfile{\s\s\s  http://www.maths.warwick.ac.uk/agt/AGTVol\thevolumenumber/agt-\thevolumenumber-\thepapernumber.abs.html}
\immediate\write\gtoutfile{\noexpand\\}
\immediate\write\gtoutfile{}
\ifx\theasciiabstract\relax
\immediate\write\gtoutfile{\theabstract}\else
\immediate\write\gtoutfile{\theasciiabstract}\fi
\immediate\write\gtoutfile{}
\immediate\write\gtoutfile{\noexpand\\}
\immediate\write\gtoutfile{}
\immediate\closeout\gtoutfile}}  

\def\maketitlepage{\makeagttitle\makeheadfile}
\let\makeshorttitle\maketitlepage
\let\maketitle\maketitlepage

\lognumber{28}
\volumenumber{1}
\volumeyear{2001}
\papernumber{28}
\pagenumbers{549}{578}
\received{18 July 2001}
\revised{5 October 2001}
\accepted{5 October 2001}
\published{14 October 2001}

\usepackage{amsmath}
\usepackage{amssymb}
\usepackage{amscd} 
\usepackage{graphics}
\newcommand{\bz}{\mathbf Z}
\newcommand{\Hom}{\operatorname{Hom}}
\newcommand{\bc}{\mathbf C}
\newcommand{\bq}{\mathbf Q}

\newcommand{\ind}{\operatorname{ind}}
\newcommand{\sign}{\operatorname{sign}}

\newcommand{\sgn}{\operatorname{sgn}} 

\newtheorem{thm}{Theorem}
\newtheorem{prop}{Proposition}
\newtheorem{cor}{Corollary}
\theoremstyle{definition}
\newtheorem{defn}{Definition}
\newtheorem{rem}{Remark}
\newtheorem*{add}{Addendum}

\newtheorem*{nota}{Notation}
\nocolon

\begin{document}
\title[The intersection forms of spin 4--manifolds]
{On the intersection forms of spin 4--manifolds\\bounded by spherical 3--manifolds}
\asciititle{On the intersection forms of spin 4-manifolds
bounded by spherical 3-manifolds}

\author{Masaaki Ue}
\address{Division of Mathematics, Faculty of Integrated Human 
Studies\\Kyoto University, Kyoto, 606-8316, Japan}
\email{ue@math.h.kyoto-u.ac.jp}
\keywords{Spin 4--manifold, spherical 3--manifold, Dirac operator }
\asciikeywords{Spin 4-manifold, spherical 3-manifold, Dirac operator }
\primaryclass{57N13}
\secondaryclass{57M50, 57M60}

\begin{abstract}
We determine the contributions of isolated singularities of spin $V$
4--manifolds to the index of the Dirac operator over them. From these
data we derive certain constraints on the intersection forms of spin
4--manifolds bounded by spherical 3--manifolds, and also on the
embeddings of the real projective planes into 4--manifolds.
\end{abstract}
\asciiabstract{We determine the contributions of isolated
singularities of spin V 4-manifolds to the index of the Dirac operator
over them. From these data we derive certain constraints on the
intersection forms of spin 4-manifolds bounded by spherical
3-manifolds, and also on the embeddings of the real projective planes
into 4-manifolds.} 
\maketitle

The 10/8--theorem \cite{Fu} and its $V$ manifold version \cite{FF} have provided several 
results about the intersection forms of spin 4--manifolds. 
For example, these theorems were used to show 
the homology cobordism invariance of 
the Neumann-Siebenmann invariant for certain Seifert homology 3--spheres 
in \cite{FFU}, and for all Seifert homology 3--spheres by 
Saveliev \cite{Sa}. 
For this purpose in \cite{FFU} we 
studied the index of the Dirac operator over spin $V$ 4--manifolds, in particular those 
with only isolated singular points whose neighborhoods are cones over lens spaces. 
The spin $V$ manifolds considered in \cite{Sa} are also of the same type, although 
they are different from those considered in \cite{FFU}. 
For a closed spin $V$ 4--manifold $X$, the index of the Dirac operator over $X$ is represented as 
\[\ind D (X)  =-(\sign X +\delta (X)) /8 ,\] 
where $\sign X$ is the signature of $X$ and $\delta (X)$ is the contribution of 
the singular points to the index of the Dirac operator, which is determined only by the data on
the neighborhoods of the singular points according to the 
V-index theorem \cite{K2}. In particular if all the singular points are isolated, $\delta (X)$ is the 
sum of the contributions $\delta (x)$ of the singular points $x$. 
In \cite{FFU} we showed that $\delta (x)$ for the case when the neighborhood of $x$ is a 
cone over a lens space is determined 
by simple recursive formulae. 
In this paper we determine the value $\delta (x)$ for 
every isolated singularity $x$, and combining such data with the 10/8 theorem, we 
derive certain information on the intersection form of a spin 4--manifold bounded by 
a spherical 3--manifold equipped with a spin structure. We also apply this to the 
embeddings of the real projective plane into 4--manifolds. 
\subsubsection*{Acknowledgments} The author thanks M. Furuta and Y. Yamada for their helpful 
comments on the earlier version of this paper (see Remarks \ref{misc}, \ref{BL}). 

\section{The $V$ manifold version of the 10/8 theorem} 
First let us recall the theorem in \cite{FF}, which coincides with the 10/8 theorem for the 
case of non-singular spin 4--manifolds. 

\begin{thm}{\rm \cite{FF}}\qua \label{10/8}
Let $Z$ be a closed spin $V$ 4--manifold with $b_1 (Z)=0$. Then either 
$\ind D(Z)=0$ or 
\[ 1-b^- (Z) \le \ind D(Z) \le b^+ (Z) -1 .\]
Since $\ind D(Z)$ is even, we have $\ind D(Z) =0$ if $b^{\pm} (Z) \le 2$. 
\end{thm}
A direct application of this theorem leads us to the following result. 

\begin{prop} \label{10/8,2}
Suppose that a 3--manifold $M$ with a spin structure $c$ bounds a spin $V$ manifold $X$ with 
isolated singularities with $b_1 (X)=0$. 

\begin{enumerate}
\item If $M$ also bounds a spin 4--manifold $Y$, then 
either $\sign (X) +\delta (X) =\sign (Y)$, or both of the following inequalities hold.
\begin{align*} 
&b^+ (Y)-9b^- (Y) \le \sign (X) +\delta (X) +8b^+ (X) -8 \\
&9b^+ (Y) -b^- (Y) \ge \sign (X) +\delta (X) -8b^- (X) +8 . 
\end{align*}
In either case, $\sign (X) +\delta (X)$ must be equal to $\sign (Y) \pmod {16}$, which is the 
Rochlin invariant $R(M,c) \pmod {16}$ of $(M,c)$ 
\item If both $b^+ (X) \le 2$ and $b^- (X) \le 2$ and 
$M$ bounds a 
$\bq$ acyclic spin 4--manifold, 
then $\sign (X) +\delta (X)=0$. 
\end{enumerate}
\end{prop}

\begin{proof}
We can assume that $b_1 (Y)=0$, for otherwise we can perform a spin surgery to get a new spin 
4--manifold $Y'$ with $b_1 (Y ')=0$, $b^{\pm} (Y)=b^{\pm} (Y ')$, and $\sign Y=\sign Y '$. 
The first claim comes from the application of Theorem \ref{10/8} to the 
index of the Dirac operator on 
$X\cup (-Y)$, given by 
\[ \ind D(X \cup (-Y))= 
-(\sign X -\sign Y +\delta (X) )/8. \] 
Note that the value in the parentheses on the right hand side must be divisible by 16 since 
the index (over $\bc$) of the Dirac operator associated with the spin structure is even. 
To prove the second claim, suppose that $M$ bounds a spin $\bq$ acyclic 4--manifold $Y$. Then 
Theorem \ref{10/8} shows that either $\sign (X) +\delta (X) = \sign (Y) =0$ or 
\[ -1 \le 1-b^- (X) \le \ind D(X\cup (-Y)) =-(\sign X +\delta (X) ) /8 \le b^+ (X) -1 \le 1 . \] 
Since $\ind D(X\cup (-Y))$ is even, we obtain the desired result. 
\end{proof}
In case of spherical 3--manifolds we obtain the following stronger result. 

\begin{prop} \label{10/8,3} 
Let $S$ be a spherical 3--manifold equipped with the spin structure $c$, and $\delta (S,c)$ be the 
contribution of the cone $cS$ over $S$ to the index of the Dirac operator $($we will show in the 
next section that such a contribution is determined only by $(S,c)$$)$. 
\begin{enumerate}
\item If $S$ bounds a spin 4--manifold $Y$, then either 
$\sign (Y)=\delta (S,c)$ or 
\[ \text{$b^+ (Y) -9b^- (Y) \le \delta (S,c) -8$ and $9b^+ (Y) -b^- (Y) \ge \delta (S,c) +8 $.} \] 
\item Suppose that for some $k$ the connected sum $kS$ of $k$ copies of $S$ $($equipped with 
the spin structure induced by $c$$)$ bounds a $\bq$ acyclic spin 4--manifold 
$($whose spin structure is an extension of the given one on $kS$$)$, then 
$\delta (S,c) =0$ $($and hence $R(S,c) \equiv 0 \pmod {16}$$)$. In particular any 
$\bz_2$ homology 3--sphere $S$ with $\delta (S) \ne 0$ $($or $R(S) \not\equiv 0 \pmod {16}$$)$ 
has infinite order in the homology cobordism group $\Theta^3_{\bz_2}$ of 
$\bz_2$ homology 3--spheres. 
\end{enumerate}
\end{prop}

\begin{proof}
Again it suffices to prove the claim for the case 
when $b_1 (Y) =0$. We can apply Proposition 1 by putting $X=cS$ to prove the first claim. 
(We will prove in the next section that $c$ extends uniquely to the spin structure on 
$cS$.) 
In this case $\sign X =b^{\pm } (X) =0$ and $\delta (X) =\delta (S, c)$. To prove the second claim 
suppose that $kS$ bounds a spin $\bq$ acyclic 4--manifold $Y$. Then applying Theorem \ref{10/8} to 
the closed spin $V$ manifold $Z$ obtained by gluing the boundary connected sum of $k$ copies of 
$cS$ and $-Y$, we have $\ind D(Z) =0$ since $b^{\pm} (Z) =0$. Since 
\[ \ind D(Z) =-(k\sign (cS) -\sign Y +k\delta (cS ))/8 \] 
and $\sign (Y)=0$, we obtain the desired result. 
\end{proof}

\section{Contributions from the cones over the spherical 
3--manifolds to the index of the Dirac operator}
Let $Z$ be a closed spin $V$ 4--manifold with a spin structure $c$ whose singularties consist of 
isolated points $\{ x_1, \dots ,x_k \}$. Then  
the $V$ index theorem \cite{K2} shows that 
the $V$ index over $\bc$ of the Dirac operator over $Z$ is described as 
\[ \ind D(Z) =\int_Z (-p_1 (Z)/24 ) +\sum_{i=1}^k \delta_D (x_i) ,\] 
where $\delta_D (x_i) $ is a contribution from the singular point $x_i$, which is described as follows. 
We omit the subscript $i$ for simplicity. Suppose that the neighborhood $N(x)$ of $x$ is represented as 
$D^4 /G$ (which is the cone over the spherical 3--manifold $S=S^3 /G$). Here $G$ is a finite subgroup 
of $SO(4)$ that acts freely on $S^3$. The restriction of $D$ to $N(x)$ is covered by a $G$ invariant 
Dirac operator $\widetilde D$ over $D^4$ and the normal bundle over $x$ in $Z$ is covered by a 
normal bundle $N$ over $0$ in $D^4$, which is identified with $\bc^2$. Then we have 
\[ \delta_D (x) =\sum_{(g) \subset (G), g\ne 1} \frac{1}{m_g} \cdot 
\frac {ch_g j^* \sigma ({\widetilde D}) }{ch_g \lambda_{-1} (N\otimes  \bc )} \] 
where $j \co \{ 0 \} \subset D^4$ is the inclusion, 
$m_g$ denotes the order of the centralizer of $g$ in $G$ and the sum on the right hand side ranges 
over all the conjugacy classes of $G$ other than the identity. On the other hand the signature of $Z$ 
(which is the index of the   
signature operator $D_{\sign}$ over $Z$) is given by 
\[ \sign (Z) =\int_Z p_1 (Z)/3 +\sum_{i=1}^k \delta_{D_{\sign}} (x_i) , \] 
where the local contribution $\delta_{D_{\sign}} (x)$ from $x$ to $\sign (Z)$ is described as 
\[ \delta_{D_{\sign}} (x) =\sum_{(g) \subset (G), g\ne 1 } \frac{1}{m_g } \cdot 
\frac{ch_g j^* \sigma({\widetilde {D}_{\sign} } )}{ch_g \lambda_{-1} (N\otimes \bc )} .\] 
Here $D_{\sign}$ over $N(x)$ is covered by a $G$ invariant signature operator $\widetilde {D}_{\sign}$ 
as before \cite{K1}. Hence we have 
\[ \ind D(Z)=-\frac 18 (\sign Z +\sum_{i=1}^k \delta (x_i) ) ,\] 
where 
\[ \delta (x) =-(\delta_{D_{\sign}} (x) +8\delta_D (x)) \] 
and we put $\delta (Z) =\sum_{i=1}^k \delta (x_i )$. 
If $N(x)$ is the cone over $S=S^3 /G$ we write $\delta (x) =\delta (S, c)$, 
where $c$ denotes the spin structure on $S$ induced from that on $Z$, since we will see later that 
$\delta (x)$ is determined completely by $(S,c)$. 
In \cite{FFU} $\delta (S, c)$ in the case when $S =L(p,q)$ is given explicitly 
as follows. The spin structure $c$ on the 
cone $cL(p,q)$ over $L(p,q)$ is determined by 
the choice of the complex line bundle $\widetilde K$ over $D^4$ that is a double covering of the 
canoncial bundle $K$ over $cL(p,q)$. Here $\widetilde K$ is the quotient space of $D^4 \times \bc$ by 
the cyclic group $\bz_{p}$ of order $p$ so that the action of the generator $g$ of $\bz_p$ is given 
by 
\[ g(z_1, z_2 ,w) =(\zeta z_1 , \zeta^q z_2 ,\epsilon \zeta^{-(q+1)/2} w) , \] 
where $\zeta =\exp (2\pi i /p)$ and 
$\epsilon =\pm 1$. There is a one-to-one correspondence between the choice of the spin structure 
on $L(p,q)$ and that of $\epsilon$. We note that every spin structure on $L(p,q)$ extends uniquely 
to that on $cL(p,q)$ and we must have $\epsilon =(-1)^{q-1} $ if $p$ is odd (see \cite{F}, \cite{FFU}). 

\begin{defn}{\rm \cite{FFU}}\qua \label{sigma}
For $L(p,q)$ with a spin structure $c$, which corresponds to the sign $\epsilon$ as above, 
$\delta (L(p,q) , c )$ equals $\sigma (q,p,\epsilon )$, which is defined by 
\begin{equation} \sigma (q, p, \epsilon )=\frac 1p \sum_{k=1}^{|p|-1} 
\left( \cot (\frac{\pi k}{p})\cot (\frac{\pi kq}{p})+2\epsilon^k 
\csc  (\frac{\pi k}{p})\csc  (\frac{\pi kq}{p}) 
\right) .
\end{equation} 
Here $p$ or $q$ may be negative under the convention $L(p,q)=L(|p| ,(\sgn p)q)$.  
\end{defn}

In \cite{FFU} we give the following characterization of $\sigma (q,p,\epsilon )$. 

\begin{prop} {\rm \cite{FFU}}\qua \label{char}
$\sigma (q,p,\epsilon )$ is an integer characterized uniquely by the 
following properties. 

\begin{enumerate}
\item $\sigma (q+cp, p, \epsilon )=\sigma (q,p, (-1)^c \epsilon )$. 
\item $\sigma (-q, p ,\epsilon )=\sigma (q,-p, \epsilon )=-\sigma (q,p,\epsilon )$.
\item $\sigma (q,1,\epsilon )=0$. 
\item $\sigma (p,q,-1 )+\sigma (q,p,-1) =-\sgn (pq)$ if $p+q \equiv 1 \pmod 2$. 
\end{enumerate}
\end{prop}

\begin{prop}{\rm\cite{FFU}}\qua \label{conti}
If $p+q \equiv 1 \pmod 2$ and $|p|>|q|$ then for a unique continued fraction 
expansion of the form 
\[ p/q =[[ \alpha_1 ,\alpha_2 .\dots ,\alpha_n ]] =
\alpha_1 -\cfrac{1}{ \alpha_2 -\cfrac {1}{\ddots -\cfrac {1}{\alpha_n }}} \]
with $\alpha_i$ even and $|\alpha_i | \ge 2$, we have 
\[ \sigma (q, p, -1) =-\sum_{i=1}^n \sgn \alpha_i .\] 
\end{prop}

\begin{cor} \label{cor}
For any coprime integers $p$, $q$ with $p$ odd and $q$ even, we have 
$\sigma (p,q, \pm 1 ) \equiv 1 \pmod 2$ and $\sigma (q,p,-1) \equiv 0 \pmod 2$. 
\end{cor}

\begin{proof}
If $|p|>|q|$ and $p$ and $q$ have opposite parity, then 
in the continued fraction expansion of 
$p/q$ in Proposition \ref{conti}
we can see inductively that $q \equiv n \pmod 2$, and hence 
$\sigma (q, p, -1) \equiv n \equiv q \pmod 2$ by Proposition \ref{conti}. 
It follows from Proposition \ref{char} 
that $\sigma (p,q,-1) \equiv p \pmod 2$. If $p$ is odd and $q$ is even then 
$\sigma (p,q,1) =\sigma (p+q, q, -1 ) \equiv p+q \equiv p \pmod 2$ also by Proposition \ref{char}. 
This proves the claim. 
\end{proof}

\begin{figure}[ht!]
\begin{center}
\includegraphics{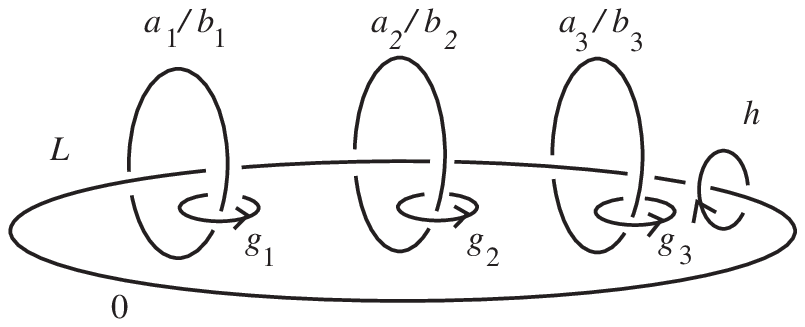}
\caption{\label{Seifert} }
\end{center}
\end{figure}

Next we consider $\delta (S, c)$ for a spherical 3--manifold $S =S^3 /G$ with 
nonabelian fundamental group $G$ with spin structure $c$. Such a manifold $S$ is a Seifert 
manifold 
over a spherical 2--orbifold $S^2 (a_1 ,a_2 ,a_3 )$ represented by the Seifert invariants of the form 
\[ S= \{ (a_1 ,b_1 ), (a_2 ,b_2 ), (a_3 ,b_3 ) \} \]
with $a_i \ge 2$, $\gcd (a_i ,b_i ) =1$ for $i=1,2,3$, $\sum_{i=1}^3 1/a_i >1$, 
$e=-\sum_{i=1}^3 b_i /a_i  \ne 0$. Here we adopt the convention in [NR] so that $S$ is 
represented by a framed link $L$ as in Figure \ref{Seifert}. 

The meridians $g_i$ and $h$ in Figure \ref{Seifert} generate $G$ with relations: 
\[  g_1^{a_1} h^{b_1} =g_2^{a_2} h^{b_2} =g_3^{a_3} h^{b_3} =g_1 g_2 g_3 =1, \quad 
[g_i ,h]= 1 \quad (i=1,2,3) . \]
The representation above is unnormalized. 
We can choose the other curves $g_i '$ homologous to $g_i +c_i h$ with 
$\sum_{i=1}^3 c_i =0$, which give an alternative representation of $S$ of the form 
\[ \{ (a_1 ,b_1 -a_1 c_1 ), (a_2, b_2 -a_2 c_2 ), (a_3 , b_3 -a_3 c_3 ) \} .\] 

Furthermore $-S$ is represented by $\{ (a_1 ,-b_1 ), (a_2 ,-b_2 ), (a_3 ,-b_3 ) \}$. 
Thus the class of the spherical 3--manifolds  
with non-abelian fundamental group up to orientation is given by the following list. 

\begin{enumerate}
\item $\{ (2,1), (2,1), (n,b) \}$ with $n\ge 2$, $\gcd (n,b) =1$, 
\item $\{ (2,1), (3,1), (3,b) \}$ with $\gcd (3,b) =1$, 
\item $\{ (2,1), (3,1), (4,b) \}$ with $\gcd (4,b) =1$, 
\item $\{ (2,1), (3,1), (5,b) \}$ with $\gcd (5,b)=1.$ 
\end{enumerate}

We also note that the above class together with the lens spaces coincides with the class of the links of the 
quotient singularities. The orientation of $S$ induced naturally by the complex 
orientation is given by choosing the signs of the Seifert invariants so that 
the rational Euler class $e$ is negative. 

\begin{defn}{\rm\cite{FFU}}\qua \label{spin}
Let $M$ be a 3--manifold represented by a framed link $L$, and let $m_i$ and $\ell_i$ be 
the meridian and the preferred longitude of the component $L_i$ of $L$ with framing 
$p_i/q_i$. Denote by $M_i$ the meridian of the newly attached solid torus along $L_i$ 
(homologous to $p_i m_i +q_i \ell_i$ in $S^3 \setminus L_i$). Then according to \cite{FFU} 
we describe a spin structure $c$ on $M$ by a homomorphism $w \in 
\Hom (H_1 (S^3 \setminus L ,\bz ),\bz_2 )$ so that 
\[ w(M_i ) := p_i w(m_i ) +q_i w(\ell_i ) +p_i q_i \pmod 2 \] 
is zero for every component $L_i$. Note that $w(m_i )=0$ if and only if $c$ extends to the 
spin structure on the meridian disk in $S^3$. 
\end{defn}

Hereafter the above homomorphism $w$ is denoted by the same symbol $c$ as the spin structure on 
$S$ if there is no danger of confusion. 
Thus the spin structures on $S=\{ (a_1, b_1 ), (a_2 ,b_2 ), (a_3 ,b_3 ) \}$ 
correspond to the elements $c\in \Hom (H_1 (S^3 \setminus L, \bz ), \bz_2 )$ satisfying 
\begin{equation} a_i c(g_i )+b_i c(h) \equiv a_i b_i \quad (i=1,2,3), \quad 
\sum_{i=1}^3 c(g_i) \equiv 0 \pmod 2 
\end{equation}

\begin{prop} \label{extend}
Every spin structure on 
the spherical 3--manifold $S$ extends uniquely to that on the cone $cS =D^4 /G$ over $S$. 
\end{prop}

\begin{proof}
The claim for a lens space was proved in \cite{F}.
We can assume that up to conjugacy $G$ is contained in $U(2) =S^3 \times S^1 /\bz_2$ (\cite{S}). 
Since the tangent frame bundle of $S$ is trivial, the associated stable $SO(4)$ bundle is reduced to 
the $U(2)$ bundle, which is represented as $S^3 \times U(2) /G$. A spin structure on $S$ corresponds to 
the double covering $S^3 \times (S^3 \times S^1 )/G \to S^3 \times U(2) /G$ for some representation 
$G \to S^3 \times S^1$ that covers the original representation of $G$ to $U(2)$. Using this representation 
we have a double covering $D^4 \times (S^3 \times S^1 )/G \to D^4 \times U(2) /G$, which gives 
a spin structure on the $V$ frame bundle over $cS=D^4 /G$. Passing to the determinant bundle (which is 
the dual of the canonical bundle of $cS$), we have a double covering of the representation 
$G \to S^1$ defined by the determinant of the element of $G$. 
Such coverings are classified by $H^1 (G, \bz_2 )=H^1 (S, \bz_2 )$. It follows that there is a one-to-one 
correspondence between the set of spin structures on $S$ and that for $cS$. This proves the claim.
\end{proof}

\begin{figure}[ht!]
\begin{center}
\includegraphics{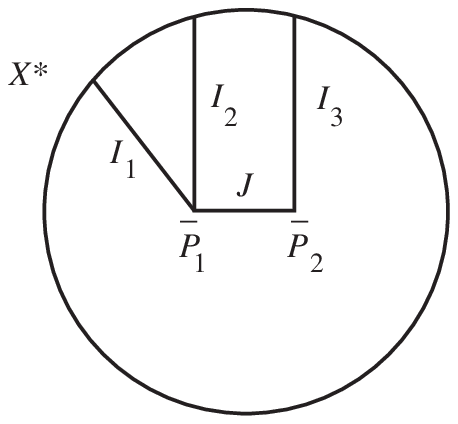}
\caption{\label{orbit} }
\end{center}
\end{figure}

Thus for a spin structure $c$ on $S$, we also denote its unique extension to $cS$ by $c$ and 
the contribution of $cS$ to the index of the Dirac operator by $\delta (S,c)$. 
To compute $\delta (S,c)$, 
we appeal to the vanishing theorem of 
the index of the Dirac operator on a certain $V$ manifold as in \cite{FFU}. 
(There is an alternative method of computing $\delta (S,c)$ by using plumbing constructions. See 
\S 3.) 
For this purpose we consider the $V$ 4--manifold $X$ with $S^1$ action and with $\partial X=S$, which 
is constructed as follows. We denote by $\pi \co X\to X^*$ the projection to the orbit space 
$X^*=X/S^1$. Suppose that $S=\{ (a_1 ,b_1 ), (a_2 ,b_2 ), (a_3 ,b_3 ) \}$ with the spin structure $c$. 
Then $X^*$ has the following properties (see Figure \ref{orbit}).  

\begin{enumerate}
\item The underlying space of $X^*$ is the 3--ball. 
\item The image of the fixed points consists of two interior points $\overline P_i$ $(i=1,2)$, and 
the image of the exceptional orbits in $X^*$ consists of three segments $I_j$ $(j=1,2,3)$ such  
that $I_j$ connects some point on $\partial X^*$ and $\overline P_1$ (for $j=1,2$) or 
$\overline P_2$ (for $j=3$). 
\item The Seifert invariant of the orbit over any point on $I_j$ except for $\overline P_i$'s is 
$(a_j ,b_j )$. 
\item The orbit over any point outside the union of $I_j$'s has a trivial stabilizer. 
\end{enumerate}

\begin{figure}[ht!]
\begin{center}\vglue -0.5in
\includegraphics{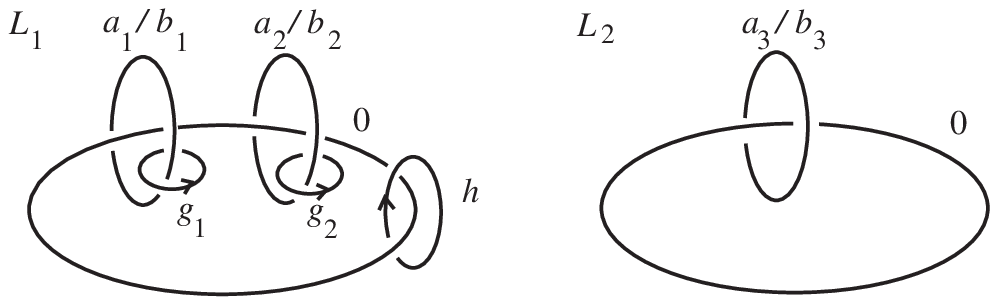} 
\caption{\label{lens} }
\end{center}
\end{figure}

Let $D_i$ be the small 4-ball neighborhood of $\overline P_i$ for $i=1,2$. Then $\pi^{-1} (D_i )$ is the cone 
over $L_i =\pi^{-1} (\partial D_i )$, where $L_i$ is represented by the framed links in Figure \ref{lens}. 
Here $L_2$ is the lens space $L(b_3, a_3 )$ represented by a $-b_3/a_3 $ surgery along the trivial knot 
with meridian corresponding to $h$, while $L_1$ is the lens space $L(Q, P)$ such that 
\begin{equation} 
\text{$Q=a_1 b_2 +a_2 b_1$, $P =a_2 v_1 +b_2 u_1$ for $u_1$, $v_1 \in \bz $ with $a_1 v_1 -b_1 u_1 =1$}. 
\end{equation}

Note that $L_1$ is represented by the 
$-Q/P$ surgery along the trivial knot whose meridian corresponds to 
\begin{equation} m =u_1 g_1 +v_1 h. \end{equation}
It follows that $X$ is a $V$ manifold with $\partial X =S$, and with two singular points 
$P_i =\pi^{-1} (\overline P_i )$ whose neighborhoods are the cones over $L_i$. 
Now according to the argument in \cite{FFU} we can check the properties of $X$. (In \cite{FFU} such a 
construction was considered when $\partial X$ is a $\bz$ homology 3--sphere. But the argument there is 
valid when $\partial X$ is a $\bq$ homology 3--sphere without any essential change.) 
Let $J$ be the segment that connects $\overline P_1$ and $\overline P_2$ in the interior of $X^*$ and 
disjoint from the interior of $I_j$ (Figure \ref{orbit}). Then $S_0 :=\pi^{-1} J$ 
is a 2--sphere and $X$ is homotopy equivalent to $S_0$. Furthermore the rational 
self intersection number of $S_0$ is given by 
\begin{equation} S_0 \cdot S_0 =a_1 a_2 /Q +a_3 /b_3 , \end{equation} 
which is nonzero if $\partial X$ is a $\bq$ homology 3--sphere. 
It follows that $b_1 (X) =0$, $b_2 (X) =1$, and $\sign X =\sgn S_0 \cdot S_0$. 
Furthermore $X$ admits a spin structure extending $c$ on $S=\partial X$ if and only if 
$c \in \Hom (H_1 (S^1 \setminus L ,\bz ),\bz_2 )$ satisfies the following conditions. 
\begin{equation} a_i c(g_i ) + b_i c(h) \equiv a_i b_i \quad (i=1,2,3), \quad \sum_{i=1}^2 c(g_i) \equiv c(g_3 )\equiv 0. 
\pmod 2 \end{equation} 
In our case $c(h)$ must be 0 since $(a_i ,b_i ) = (2,1)$ for some $i$. 
We will see later that we can arrange the Seifert invariants for any given $(S,c)$ so that they 
satisfy these conditions. If we put  
\[ \widehat X =X \cup (-cS ) ,\] 
then by Proposition \ref{extend} $\widehat X$
is a closed spin $V$ 4--manifold  with $b_1 (\widehat X )=0$, $b_2 (\widehat X )=1$ and 
$\sign \widehat X =\sgn S_0 \cdot S_0 $, such that $\widehat X$ has (at most) three singular points 
whose neighborhoods are the cones over $L_1$, $L_2$, and $-S$. 
Here the Seifert invariants for $-S$ are given by 
$\{ (a_1 ,-b_1 ), (b_2 ,-b_2 ), (a_3 ,-b_3 ) \}$ with respect to the curves $g_i$ and $-h$, 
and we can consider the spin structure on $-S$ induced from $c$, which is given by the 
same homomorphism in 
$\Hom (H_1 (S^3 \setminus L, \bz ),\bz_2 )$ as $c$ 
and is denoted by $-c$. Then the spin structure on $-S$ induced from $\widehat X$ 
is $-c$. Moreover 
the spin structures on $L_1$ and $L_2$ 
induced from 
that on $\widehat X$ correspond to $c(h)$ and $c(m)$ respectively, where 
\begin{equation} c(m) \equiv u_1 c(g_1 )+v_1 c(h) +u_1 v_1 \pmod 2 \end{equation} 
(see \cite{FFU}). 
Then the argument in \cite{FFU}, Proposition 3 shows that 
\begin{equation} 
\delta (L_1 ,c)=\sigma (P, Q, (-1)^{(c(m) -1) }), \quad \delta (L_2 ,c)=\sigma (a_3 ,b_3 ,(-1)^{(c(h)-1)}).
\end{equation}

Thus from Theorem 1 \cite{FF} we deduce 
\begin{equation*}
\begin{split} 
0&=\ind D(\widehat X) \\
&=
-(\sign \widehat X +\sigma (P,Q, (-1)^{(c(m)-1)} )+\sigma (a_3, b_3 ,(-1)^{(c(h)-1)} ) +
\delta (-S ,-c))/8 , 
\end{split}
\end{equation*}

Thus we can see that 
\begin{equation} 
\delta (-S, -c )=-\delta (S,c)
\end{equation}

and hence (using the fact that $c(h)=0$), 
\begin{equation} 
\delta (S,c)=\sgn S_0 \cdot S_0 +\sigma (P,Q, (-1)^{(c(m)-1)} )+\sigma (a_3 ,b_3 ,-1 ). 
\end{equation} 

Now we apply this result to compute $\delta (S,c)$ by constructing the above $X$ 
associated with the Seifert invariants of $S$, which is rearranged if necessary. We denote 
by $\{ (a_1, b_1 ), (a_2 ,b_2 ), (a_3 ,b_3 ) \}$ the rearraged Seifert invariants and 
the corresponding meridian curves in Figure \ref{Seifert} by $g_i '$ ($h$ remains unchanged). 
Hereafter we write 
the data of the required $X$ by giving the rearranged Seifert invariants, the values of 
$Q$, $P$, $m$, and $S_0 \cdot S_0$. 

\subsection{Case 1\qua $S=\{ (2,1), (2,1), (n,b ) \} $}
\subsubsection{ $n$ is odd and $b$ is even} 
In this case the spin structure $c$ on $S$ satisfies 
\begin{equation} 
c(h) \equiv c(g_3 ) \equiv 0, \quad c(g_1) \equiv c(g_2 ) =\epsilon \pmod 2  
\end{equation}
where 
$\epsilon$ is arbitrary. Then $X$ associated with $(a_1 ,b_1 )=(a_2 ,b_2 )=(2,1)$, $(a_3 ,b_3 )=
(n,b)$, 
$Q=4$, $P=3$, $m=g_1 '+h$, $S_0 \cdot S_0 =(n+b)/b$ shows that 
\[ \delta (S,c)=\sgn (n+b)b +\sigma (3,4,(-1)^{\epsilon} ) +\sigma (n,b,-1) .\] 

Since $\sigma (3,4,1)=\sigma (7,4,-1) =-\sigma (4,7,-1) -1$ and 
$4/3 =[[ 2,2,2]]$, $7/4 =[[2,4]]$, 
we deduce from Propositions \ref{char} and \ref{conti} that 
\begin{equation} \label{delta value} 
\delta (S,c) =
\begin{cases}
\sigma (n,b,-1) &(\epsilon =0, \quad -n<b<0) \\
\sigma (n,b,-1 )+2 &(\epsilon =0, \quad \text{$b>0$ or $b<-n$}) \\
\sigma (n,b,-1)-4 &(\epsilon =1, \quad -n<b<0 ) \\
\sigma (n,b,-1)-2 &(\epsilon =1, \quad \text{$b>0$ or $b<-n$}). 
\end{cases}
\end{equation}

In either case $\delta (S,c)$ is odd by Corollary \ref{cor}. 
\subsubsection{$n$ and $b$ are odd} 
In this case $c$ is given by 
\begin{equation} c(h)=0,\quad  c(g_3) \equiv c(g_1 )+c(g_2 ) \equiv 1 \pmod 2 .
\end{equation} 
It suffices to 
consider the case when $c(g_1) \equiv 1$ and $c(g_2 )\equiv 0 \pmod 2$, since 
we have a self-diffeomorphism of $S$ mapping $(g_1, g_2 ,g_3 ,h)$ to 
$(-g_2, -g_1 ,-g_3 ,-h )$. 
Thus $X$ associated with $(a_1 ,b_1 )=(a_3 ,b_3 ) =(2,1)$ and $(a_2 ,b_2 )=(n,b)$, 
$Q=n+2b$, $P=n+b$, $m=g_1 ' +h$ and $S_0 \cdot S_0 =4(n+b)/(n+2b)$ shows that ($X$ has 
only one singular point since $L(1,2)$ is the 3--sphere) 
\begin{equation}
\begin{split} 
\delta (S, c)&=\sgn (n+b)(n+2b) +\sigma (n+b, n+2b ,-1)  \\ 
&=-\sigma (n+2b,n+b, -1)=-\sigma (-n ,n+b, -1)=\sigma (n,n+b, -1) 
\end{split} 
\end{equation}
Again $\delta (S,c)$ is odd in this case by Corollary \ref{cor}. 

\subsubsection{$n$ is even} 
In this case $c$ satisfies 
\begin{equation} c(h)\equiv 0, \quad c(g_1 )+c(g_2 )+c(g_3 ) \equiv 
0 \pmod 2. 
\end{equation}
Since at least one of $c(g_i)$ is zero and there is a self-diffeomorphism of $S$ 
exchanging $g_1$ and $g_2$ up to orientation as before, it suffices to consider the following subcases. 

\noindent
$(i)$ $c(h)\equiv c(g_3 ) \equiv 0$, \quad $c(g_1 )=c(g_2 ) =\epsilon \pmod 2$. 

Then $X$ associated with $(a_1 ,b_1 )=(a_2 ,b_2 )=(2,1)$ and $(a_3 ,b_3 )=(n,b)$, 
$Q=4$, $P=3$, $m=g_1 '+h$, and $S_0 \cdot S_0 =(n+b)/b$ shows that 
\[ \delta (S ,c)=\sgn (n+b)b +\sigma (3,4,(-1)^{\epsilon } )+\sigma (n,b,-1) \] 

Hence $\delta (S,c)$ is represented by the same equation as in Case (2.1.1) \ref{delta value}.

$(ii)$ $c(h)\equiv c(g_2 ) \equiv 0$, \quad $c(g_1 )\equiv c(g_3 )\equiv 1 \pmod 2$. 

If $Q=n+2b \ne 0$ (i.e., if $(n,b) \ne (2,-1)$) then 
$X$ associated with $(a_1 ,b_1 )=(a_3 ,b_3 )=(2,1)$ and $(a_2 ,b_2 )=(n,b)$, 
$Q=n+2b$, $P=n+b$, $m=g_1 '+h$, $S_0 \cdot S_0 = 
4(n+b)/(n+2b)$ 
shows that 
\[ \delta (S,c)=\sgn (n+b)(n+2b) +\sigma (n+b ,n+2b, -1 ) .\] 
By Proposition \ref{char} the right hand side equals  
\[
-\sigma (n+2b, n+b ,-1) 
=-\sigma (-n ,n+b, -1) =\sigma (n,n+b, -1) .
\]

For the case when $(n,b)=(2,-1)$, we consider another 
representation of $S$ of the form 
$\{ (2,1), (2,-3), (2,3) \}$ 
with respect to the curves $g_1 '=g_1$, $g_2 '=g_2 +2h$, and $g_3 '=g_3 -2h $. 
Since $c(g_1 ') \equiv c(g_1 ) \pmod 2$ and 
\[ c(g_2 ')\equiv c(g_2 )+2c(h) +2 \equiv c(g_2 ), \quad 
c(g_3 ')\equiv c(g_3 )-2c(h)-2 \equiv c(g_3 ) ,\] 
considering $X$ with 
$(a_1, b_1 )=(2,1)$, $(a_2 ,b_2 )=(2,3)$, $(a_3 ,b_3 )=(2,-3)$, 
$Q=8$, $P=5$, $m=g_1 ' +h$, and $S_0 \cdot S_0 =-1/6$, we have 
\begin{equation*} 
\delta (S,c)=-1 +\sigma (5,8,-1 )+\sigma (2,-3,-1) =0 .
\end{equation*} 

It follows that in either case 
\begin{equation} \delta (S,c) = 
\sigma (n,n+b, -1) .
\end{equation} 

We also note that there are some overlaps in the above list if 
we also consider $(-S, -c)$. In fact we have 
\begin{equation} \label{overlap 1}
\{ (2,-1), (2,-1), (n,-b) \} =
\{ (2,1), (2,1), (n, -2n-b) \} .
\end{equation}

\subsection{Case 2\qua $S=\{ (2,1), (3,1), (3,b) \} $}
In this case $S$ is a $\bz_2$ homology 3--sphere and $c$ is uniquely determined. 
\subsubsection{$b$ is even}
In this case $c$ satisfies $c(h)\equiv c(g_3 )\equiv 0$, $c(g_1 )=c(g_2 )\equiv 1
\pmod 2$. Thus $X$ associated with $(a_1 ,b_1 )=(2,1)$, $(a_2 ,b_2 )=(3,1)$, $(a_3 ,b_3 )=
(3,b)$, $Q=5$, $P=4$, $m=g_1 '+h$, $S_0 \cdot S_0 =
(6b+15)/5b$ shows that 
\[ \delta (S, c)=\sgn (2b+5)b +\sigma (4,5,-1) +\sigma (3,b,-1) .\] 

Here we must have $b=6k \pm 2$ for some $k$, and if $k\ne 0$, 
\[ 
(6k+2)/3 =[[2k, -2, -2]], \quad (6k-2)/3 =[[2k ,2,2 ]] 
\] 
and hence 
\begin{equation} \delta (S,c) =
\begin{cases} 
-\sgn b -1 &(\text{$b=6k+2$ for some $k$}) ,\\
-\sgn b-5 &(\text{$b=6k-2$ for some $k\ne 0$}), \\
-6 &(b=-2).  
\end{cases}
\end{equation}

\subsubsection{$b$ is odd} 
Consider a representation of $S$ of the form 
$\{ (2,-1), (3,1), (3, b+3 ) \}$.  
Then we have  
$c(h )\equiv c(g_3 ') \equiv 0$, $c(g_1 ')\equiv c(g_2 ') \equiv 1 
\pmod 2$. Thus $X$ associated with $(a_1 ,b_1 )=(2,-1)$, $(a_2 ,b_2 )=(3,1 )$, $(a_3 ,b_3 )
=(3, b+3 )$, $Q=-1$, $P=2$, $m=-g_1 ' +h$, $S_0 \cdot S_0 =
(3-6(b+3))/(b+3)$ shows that 
\[ \delta (S,c)=-\sgn (2b+5)(b+3) +\sigma (3,b+3 ,-1) .\] 
Here we must have $b=6k\pm 1$ for some $k$. Since 
\[ 
(6k+4)/3=[[2(k+1), 2,2]] \ \ (k\ne -1), \quad 
(6k+2)/3=[[2k, -2,-2]] \ \  (k\ne 0) , 
\] 
we can see that 
\begin{equation} \delta (S,c) =
\begin{cases}
-\sgn b -3 &(\text{$b=6k+1$ for some $k$}), \\
-\sgn b+1 &(\text{$b=6k-1$ for some $k\ne 0$}), \\
0 &(b=-1) .
\end{cases}
\end{equation}

We also note that 
\begin{equation} \label{overlap 2}
\{ (2,-1), (3,-1), (3,-6k-2) \} 
=\{ (2,1), (3,1), (3,-6(k+1)-1 ) \}. 
\end{equation}

\subsection{$S =\{ (2,1),(3,1),(4,b) \} $}
In this case $c$ satisfies 
\begin{equation}
c(h)\equiv 0, \quad c(g_2 )\equiv 1, \quad 
c(g_1 ) +c(g_3 ) \equiv 1 \pmod 2. 
\end{equation}
\subsubsection{$c(g_1 ) \equiv c(g_2 ) \equiv 1$, $c(g_3 )\equiv c(h) \equiv 0 \pmod 2$}
Considering $X$ with $(a_1 ,b_1 )=(2,1)$, $(a_2 ,b_2 )=(3,1)$, $(a_3 ,b_3 )=(4,b)$, $Q=5$, $P=4$, 
$m=g_1 '+h$, $S_0 \cdot S_0 = (6b+20)/(5b)$, we have 
\[ \delta (S,c)=\sgn (3b+10)b +\sigma (4,5,-1)+\sigma (4,b, -1) .\] 
Here we must have $b=8k \pm 1$ or $b=8k\pm 3$ for some $k$. Since 
\begin{align*}  
&(8k+1)/4=[[2k,-4]] , \quad 
(8k-1)/4=[[2k,4]], \\
&(8k+3)/4=[[2k,-2,-2,-2]], \quad (8k-3)/4=[[2k,2,2,2]] 
\end{align*}
for $k\ne 0$, we can see that 
\begin{equation}
\delta (S,c) =
\begin{cases}
-\sgn b-2 &\text{($b=8k+1$ for some $k$)}, \\
-\sgn b-4 &\text{($b=8k-1$ for some $k \ne 0$)}, \\
-\sgn b &\text{($b=8k+3$ for some $k$)}, \\
-\sgn b-6 &\text{($b=8k-3$ for some $k\ne 0$)}, \\
-5 &(b=-1) ,\\
-7 &(b=-3) .
\end{cases}
\end{equation} 

\subsubsection{$c(g_1 )\equiv c(h) \equiv 0$, $c(g_2 )\equiv c(g_3 ) \equiv 1 \pmod 2$}
In this case $X$ associated with $(a_1 ,b_1 )=(3,1)$, $(a_2 ,b_2 )=(4,b)$, $(a_3 ,b_3 )=(2,1)$, 
$Q=3b+4$, $P=2b+4$, $m=2g_1 '+h$, $S_0 \cdot S_0 =
(6b+20)/(3b+4)$ shows that 
\[ \delta (S,c)=\sgn (3b+4)(3b+10) +\sigma (2b+4, 3b+4 ,-1) . \] 
Here $(3b+4)/(2b+4)$ equals 
\begin{alignat*}{2} 
&[[2,2,2(k+1), 2,2,2]] &\quad &\text{(if $b=8k+1$ for $k\ne -1$)}, \\
&[[2,2,2k ,-2,-2,-2]] &\quad &\text{(if $b=8k-1$ for some $k\ne 0$)}, \\
&[[2,2,2(k+1),4 ]] &\quad &\text{(if $b=8k+3$ for $k\ne -1$)}, \\
&[[2,2,2k, -4]] &\quad &\text{(if $b=8k-3$ for some $k\ne 0$)}, \\
&[[2,4,2,2]] &\quad &(\text{if $b=-7$}), \\
&[[2,6]] &\quad &(\text{if $b=-5$}), \\
&[[2,-2]] &\quad &(\text{if $b=-3$}) .
\end{alignat*}

Hence we have 
\begin{equation}
\delta (S,c) =
\begin{cases}
-\sgn b-4 &(b=8k+1), \\
-\sgn b+2 &\text{($b=8k-1$ for $k\ne 0$)}, \\
-\sgn b-2 &(b=8k+3), \\
-\sgn b &\text{($b=8k-3$ for $k\ne 0$ )}, \\
-1 &(b=-3) ,\\
1 &(b=-1). 
\end{cases}
\end{equation}

\subsection {$S=\{ (2,1), (3,1), (5,b) \} $ }
In this case the spin structure on $S$ is unique. 
\subsubsection{$b$ is even}
In this case $c$ satisfies 
$c(h) \equiv c(g_3 )\equiv 0$, 
$c(g_1 )\equiv c(g_2 )\equiv 1 \pmod 2$. Then $X$ associated with 
$(a_1 ,b_1 )=(2,1)$, $(a_2 ,b_2 )=(3,1)$, $(a_3 ,b_3 )=(5,b)$, 
$Q=5$, $P=4$, $m=g_1 '+h$, $S_0 \cdot S_0 =(6b+25)/(5b)$ shows that 
\[ \delta (S,c) =\sgn (6b+25)b +\sigma (4,5,-1) +\sigma (5,b,-1). \] 

Here we must have $b=10k\pm2$, or $10k\pm4$ for some $k$. Since for $k\ne 0$ 
\begin{align*}
&(10k+2)/5 =[[2k,-2,2]] ,  \quad (10k-2)/5=[[2k,2,-2]], \\
&(10k+4)/5=[[2k,-2,-2,-2,-2]] ,\quad (10k-4)/5=[[2k,2,2,2,2]] ,
\end{align*}

we have 
\begin{equation}
\delta (S,c) =
\begin{cases}
-\sgn b-3 &\text{($b=10k \pm 2$ other than $-2$)}, \\
-\sgn b+1 &(b=10k+4), \\
-\sgn b -7 &\text{($b=10k-4$ other than $-4$)}, \\
-4 &(b=-2), \\
-8 &(b=-4). 
\end{cases}
\end{equation}

\subsubsection{$b$ is odd}
Consider the Seifert invariants of $S$ of the form 
$\{ (2,-1), (3,1), (5,b+5) \}$.
Then $c(h )\equiv c(g_3 ') \equiv 0$, $c(g_1 ')\equiv c(g_2 ') \equiv 1 \pmod 2$. 
Hence $X$ associated with $(a_1, b_1 )=(2,-1)$, $(a_2 ,b_2 )=(3,1)$, $(a_3 ,b_3 )=(5, b+5 )$, 
$Q=-1$, $P=2$, $m=-g_1 ' +h$, $S_0 \cdot S_0 =(5-6(b+5))/(b+5)$ shows that 
\[ \delta (S,c)=-\sgn (b+5)(6b+25) +\sigma (5,b+5 ,-1) .\] 

Here we must have $b=10k\pm 1$ or $10k\pm 3$ for some $k$. Since 
\begin{align*}
&(10k+6)/5=[[2(k+1),2,2,2,2]] \quad (k\ne -1) ,\\
&(10k+4)/5=[[2k,-2,-2,-2,-2]] \quad (k\ne 0) ,\\
&(10k+8)/5=[[2(k+1),2,-2]] \quad (k\ne -1), \\
&(10k+2)/5=[[2k,-2,2]] \quad (k \ne 0) ,
\end{align*}
we have 
\begin{equation}
\delta (S,c) =
\begin{cases}
-\sgn b-5 &\text{($b=10k+1$ for some $k$)} \\
-\sgn b +3 &\text{($b=10k-1$ for some $k\ne 0$)} \\
-\sgn b-1 &\text{($b=10k\pm 3$ for some $k$ and $b \ne -3$)}, \\
2 &(b=-1), \\
-2 &(b=-3) . 
\end{cases}
\end{equation}
 
Now we remove the overlaps (\ref{overlap 1}, \ref{overlap 2}) from the above results by giving the data 
only for the spherical 3--manifolds with negative rational Euler class.

\begin{prop} \label{delta}
The value $\delta =\delta (S,c)$ for a spherical 3--manifold $S$ with negative rational Euler 
class and its spin structure $c$ 
is given by the following list. Note that $\delta (-S, -c) =-\delta (S,c)$. 
Except for the lens spaces, we give the list of the 
Seifert invariants for $S$, 
the set of the values $(c(g_1 ), c(g_2 ), c(g_3 ))$ of $c$, and $\delta$. Here 
$g_i$ and $h$ are the meridians of the framed link associated with the Seifert invariants as in 
Figure \ref{Seifert}. 
We omit $c(h)$ since it is always zero. In the list below, the data of $c$ is omitted when 
$S$ is a $\bz_2$ homology sphere $($cases $(3)$ and $(5)$$)$, and $\epsilon$ is $\pm1$. 
 
\begin{flushleft} 
{$(1)$ $S=L(p,q)$ with $p>q>0$.} 
\end{flushleft}
In this case $\delta (L(p,q) ,c) =\sigma (q,p,\epsilon )$ where 
the relation between $c$ and $\epsilon$ is explained in the paragraph before Definition \ref{sigma}. 
We also note that if $L(p,q)$ is represented by the $-p/q$--surgery along the trivial knot $O$, 
then the spin structure given by $c\in \Hom ( H_1 (S^3 \setminus O, \bz ), \bz_2 )$ explained 
as in Definition \ref{spin} satisfies $\epsilon \equiv c(\mu) -1 \pmod 2$ with 
respect to the above correspondence, where $\mu$ is the 
meridian of $O$ $($see \cite{FFU}$)$.  

\begin{flushleft}
{$(2)$ $S=\{ (2,1), (2,1), (n,b) \}$.} 
\end{flushleft}
\begin{tabular}{llll}
 & $S$ & $c$ & $\delta$ \\
\hline
$(2-1)$ & $n$ odd, $b$ even, $-n<b<0$ & $(0,0,0)$ & $\sigma (n,b,-1)$ \\
$(2-2)$ & $n$ odd, $b$ even, $-n<b<0$ & $(1,1,0)$ & $\sigma (n,b,-1)-4$ \\
$(2-3)$ & $n$ odd, $b$ even, $b>0$ & $(0,0,0)$ & $\sigma (n,b,-1)+2$ \\
$(2-4)$ & $n$ odd, $b$ even, $b>0$ & $(1,1,0)$ & $\sigma (n,b,-1)-2$ \\
$(2-5)$ & $n$, $b$ odd, $n+b>0$ & $(\epsilon, 1-\epsilon ,1 )$ & $\sigma (n,n+b ,-1)$ \\
$(2-6)$ & $n$ even, $-n<b<0$ & $(0,0,0)$ & $\sigma (n,b,-1)$ \\
$(2-7)$ & $n$ even, $-n<b<0$ & $(1,1,0)$ & $\sigma (n,b,-1)-4$ \\
$(2-8)$ & $n$ even, $b>0$ & $(0,0,0)$ & $\sigma (n,b,-1) +2$ \\
$(2-9)$ & $n$ even, $b>0$ & $(1,1,0)$ & $\sigma (n,b,-1) -2$ \\
$(2-10)$ & $n$ even, $n+b >0$ & $(\epsilon ,1-\epsilon, 1)$ & $\sigma (n,n+b, -1)$ 
\end{tabular}

\begin{flushleft}
{$(3)$ $S$ is a Seifert fibration over $S^2 (2,3,3)$.} 
\end{flushleft}
\begin{tabular}{llr} 
 & $S$ & $\delta$ \\
 \hline
$(3-1)$ & $\{ (2,1),(3,1),(3,6k+2) \}$, $k\ge 0$ & $-2$ \\ 
$(3-2)$ & $\{ (2,-1), (3,-1), (3,-6k-2 )\}$, $k\le -1$ &  $0$ \\
$(3-3)$ & $\{ (2,1), (3,1), (3,6k-2) \}$, $k\ge 0$ & $-6$ \\ 
$(3-4)$ & $\{ (2,-1), (3,-1), (3, -6k+2 ) \}$, $k <0$ & $4$ \\ 
$(3-5)$ & $\{ (2,1), (3,1), (3,6k+1 ) \}$, $k \ge 0$ & $-4$ \\ 
$(3-6)$ & $\{ (2,-1), (3,-1), (3,-6k-1) \} $, $k <0$ & $2$ \\
\end{tabular} 

\begin{flushleft}
{$(4)$ $S$ is a Seifert fibration over $S^2 (2,3,4)$.}
\end{flushleft}
\begin{tabular}{lllr} 
 & $S$ & $c$ & $\delta$ \\
 \hline
$(4-1)$ & $\{ (2,1),(3,1),(4,8k+1)\}$, $k\ge 0$ & $(1,1,0)$ & $-3$ \\
$(4-2)$ & $\{ (2,1), (3,1), (4,8k+1) \}$, $k\ge 0$ & $(0,1,1)$ & $-5$ \\
$(4-3)$ & $\{ (2,-1),(3,-1),(4,-8k-1) \}$, $k <0$ & $(1,1,0)$ & $1$ \\
$(4-4)$ & $\{ (2,-1), (3,-1), (4,-8k-1) \}$, $k<0$ & $(0,1,1)$ & $3$ \\
$(4-5)$ & $\{ (2,1), (3,1), (4,8k-1) \}$, $k\ge 0$ & $(1,1,0)$ & $-5$ \\
$(4-6)$ & $\{ (2,1),(3,1),(4,8k-1) \}$, $k\ge 0$ & $(0,1,1)$ & $1$ \\
$(4-7)$ & $\{ (2,-1), (3,-1), (4,-8k+1) \}$, $k<0$ & $(1,1,0)$ & $3$ \\
$(4-8)$ & $\{ (2,-1), (3,-1), (4,-8k+1) \}$, $k<0$ & $(0,1,1)$ & $-3$ \\
$(4-9)$ & $\{ (2,1), (3,1), (4,8k+3) \}$, $k\ge 0$ & $(1,1,0)$ & $-1$ \\
$(4-10)$ & $\{ (2,1), (3,1), (4,8k+3) \}$, $k\ge 0$ & $(0,1,1)$ & $-3$ \\
$(4-11)$ & $\{ (2,-1), (3,-1), (4,-8k-3) \}$, $k<0$ & $(1,1,0)$ & $-1$ \\
$(4-12)$ & $\{ (2,-1), (3,-1), (4,-8k-3) \}$, $k<0$ & $(0,1,1)$ & $1$ \\
$(4-13)$ & $\{ (2,1), (3,1), (4,8k-3) \}$, $k\ge 0$ & $(1,1,0)$ & $-7$ \\
$(4-14)$ & $\{ (2,1), (3,1), (4,8k-3) \}$, $k\ge 0$ & $(0,1,1)$ & $-1$ \\
$(4-15)$ & $\{ (2,-1), (3,-1), (4,-8k+3) \}$, $k<0$ & $(1,1,0)$ & $5$  \\
$(4-16)$ & $\{ (2,-1), (3,-1), (4,-8k+3) \}$, $k<0$ & $(0,1,1)$ & $-1$ 
\end{tabular} 

\begin{flushleft}
{$(5)$ $S$ is a Seifert fibration over $S^2 (2,3,5)$.} 
\end{flushleft}
\begin{tabular}{llr}
 & $S$ & $\delta$ \\
\hline 
$(5-1-\epsilon )$ & $\{ (2,1), (3,1), (5,10k +2\epsilon ) \}$, $k\ge 0$ & $-4$ \\
$(5-2-\epsilon )$ & $\{ (2,-1), (3,-1), (5,-10k - 2\epsilon ) \}$, $k<0$ & $2$ \\
$(5-3)$ & $\{ (2,1), (3,1), (5, 10k+4 ) \}$, $k\ge 0$ &  $0$ \\
$(5-4)$ & $\{ (2,-1), (3,-1), (5,-10k-4) \}$, $k<0$ & $-2$ \\
$(5-5)$ & $\{ (2,1), (3,1), (5,10k-4) \}$, $k\ge 0$  & $-8$ \\
$(5-6)$ & $\{ (2,-1), (3,-1), (5,-10k+4) \}$, $k<0$ & $6$ \\
$(5-7)$ & $\{ (2,1), (3,1), (5,10k+1) \}$, $k\ge 0$ & $-6$ \\
$(5-8)$ & $\{ (2,-1), (3,-1), (5,-10k-1) \}$, $k<0$ & $4$ \\
$(5-9)$ & $\{ (2,1), (3,1), (5,10k-1) \}$, $k\ge 0$ & $2$ \\
$(5-10)$ & $\{ (2,-1), (3,-1), (5,-10k+1) \}$, $k<0$ & $-4$ \\
$(5-11-\epsilon )$ & $\{ (2,1), (3,1), (5,10k + 3\epsilon ) \}$, $k \ge 0$ & $-2$ \\
$(5-12-\epsilon )$ & $\{ (2,-1), (3,-1), (5,-10k - 3\epsilon ) \}$, $k<0$ & $0$ 
\end{tabular}
\end{prop}

\section{Some applications}

\begin{figure}[ht!]
\begin{center}
\includegraphics{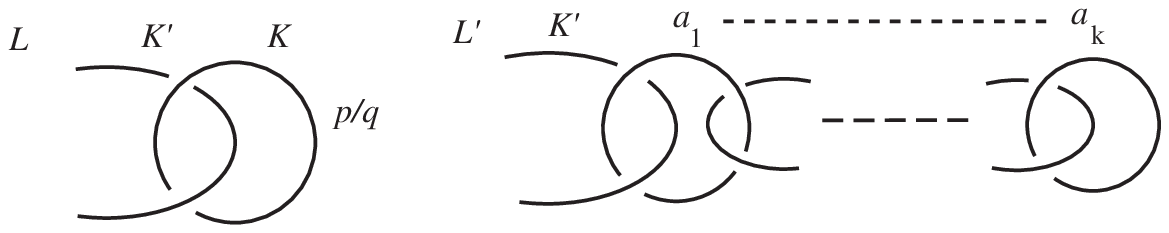}
\caption{\label{framed} }
\end{center}
\end{figure}

Let us start with some (well-known) results for later use. 
\begin{prop} \label{expand}
\begin{enumerate}
\item 
Suppose that a spin 4--manifold $Y$ is represented by a framed link $L$ with even framings. Then 
the spin structure on $\partial Y$ is induced from that on $Y$ if and only if it is represented by 
the zero homomorphism of $\Hom (H_1 (S^3 \setminus L, \bz ), \bz_2 )$. 
\item Let $M$ be a 3--manifold represented by a framed link $L$ in Figure \ref{framed}, whose framing for the 
component $K$ is given by $p/q$ for coprime $p$, $q$ with opposite parity. Suppose that 
a spin structure $c$ on $M$ is 
represented by $c\in \Hom (H_1 (S^3 \setminus L, \bz ), \bz_2 )$ with $c(\mu ) =c(\mu ')=0$ for 
meridians $\mu$ of $K$ and $\mu '$ of $K'$. 
Then the 3--manifold $M'$ represented by a framed link $L'$ in Figure \ref{framed}, where 
$p/q =[[a_1 ,\dots ,a_k ]]$ for even $a_i$, $a_i \ne 0$ is diffeomorphic to $M$, so that 
$c$ corresponds to $c' \in \Hom (H_1 (S^3 \setminus L' ,\bz ), \bz_2 )$ with 
$c' (\mu_i ) =0$ for any meridian $\mu _i$ of the new components of framing $a_i$, and $c ' (\mu '') 
=c(\mu '' )$ for a meridian $\mu ''$ of any common component of $L$ and $L'$. 
\end{enumerate}
\end{prop}

\begin{proof}
If the spin structure $c$ on $\partial Y$ is induced from that on $Y$, then the 
associated element of $\Hom (H_1 (S^3 \setminus L ,\bz ), \bz_2 )$ is 
zero since $c$ extends to that on $S^3$. Conversely if $c$ is zero, $c$ extends to the spin structure on 
$S^3$, and hence on the 4--ball, while there is no obstruction to extending $c$ to that on the 2--handles 
attached to the 4-ball since all the framings are even. 
This proves the first claim. To see the second claim note that there is a diffeomorphism between 
$M$ and $M'$ such that $\mu$ and $\mu '$ correspond to the meridians $\mu_i $ by the following relations. 
\[ 
\begin{pmatrix} -a_i & -1 \\ 1 & 0 \end{pmatrix} 
\begin{pmatrix} \mu_i \\ \mu_{i-1} \end{pmatrix} = 
\begin{pmatrix} \mu_{i+1} \\ \mu_i \end{pmatrix} \ (i\le k-1), \quad 
\begin{pmatrix} -a_k & -1 \\ 1 & 0 \end{pmatrix} 
\begin{pmatrix} \mu_k \\ \mu_{k-1} \end{pmatrix} = 
\begin{pmatrix} \widetilde{\mu} \\ \widetilde{\lambda} \end{pmatrix} 
\] 
where $(\mu_1 ,\mu_0 )=(\mu ,\mu ')$ and $(\widetilde{\mu}, \widetilde{\lambda} )$ is a pair of 
a meridian and a longitude for a newly attached solid torus along $K$.  
Since all $a_i$ are even, we have $c' (\mu_i ) =0$ for every $i$. 
\end{proof}

Next we consider plumbed 4--manifolds bounded by spherical 3--manifolds. Let $P(\Gamma )$ be a 
plumbed 4--manifold associated with a weighted tree graph $\Gamma$. Let $x_i$ be the 
generators of $H_2 (P(\Gamma ), \bz )$ corresponding to the vertices $v_i$ 
$(i\in I)$ of $\Gamma$, and $\mu_i $ be the meridian of the component associated with $v_i$ of a 
framed link $L_{\Gamma}$ naturally corresponding to $P(\Gamma )$. For every spin structure $c$ on 
$\partial P(\Gamma )$, there exists a Wu class 
$w$ of $P(\Gamma )$ associated with $c$ of the form 
$w=\sum_{i\in I} \epsilon_i x_i $ with $\epsilon_i =0$ or $1$ such that 
\[ w \cdot x_i \equiv x_i \cdot x_i \pmod 2 \quad (i\in I ) \] 
where $c$ corresponds to an element $w \in \Hom (H_1 (S^3 \setminus L_{\Gamma} ,\bz ), \bz_2 )$ 
(we use the same symbol $w$ since there is no danger of confusion) 
satisfying $w(\mu_i ) =\epsilon_i $. The set of $v_i$ with $\epsilon_i =1$ in the above representation 
of $w$ is called the Wu set (\cite{Sa}). It is well known that no adjacent vertices in $\Gamma$ both belong 
to the same Wu set. Moreover the spin structure $c$ extends to that on the complement 
in $P(\Gamma )$ of the union of 
$P(v_i )$ for $v_i$ in the Wu set. The following proposition is a generalization of the result for 
lens spaces in \cite{Sa}. 

\begin{prop} \label{wu}
Suppose that a spherical 3--manifold $S$ bounds a 
plumbed 4--manifold $P(\Gamma )$. For any spin structure 
$c$ on $S$, we have $\delta (S,c) =\sign P(\Gamma ) -w\cdot w$ for the associated Wu class $w \in H_2 (P(\Gamma ), \bz )$. 
In particular if $P(\Gamma )$ is spin and $c$ is the spin structure inherited from that on $P(\Gamma )$, 
we have $\delta (S,c )=\sign P(\Gamma )$. 
\end{prop}

\begin{proof}
It suffices to consider the case when $\Gamma$ is reduced, for otherwise by blowing down processes 
we obtain a reduced graph $\Gamma '$ such that $S=\partial P(\Gamma )=\partial P(\Gamma ' )$ and 
the Wu class $w'$ of $P(\Gamma ' )$ associated with $c$ satisfies $\sign P(\Gamma )-w\cdot w =
\sign P(\Gamma ') -w' \cdot w'$. 
In the case of lens spaces, this claim follows from the result in \cite{Sa} under the correspondence 
of $\sigma (q,p, \pm 1)$ and $\delta (L(p,q), c)$. If $S$ is not a lens space, 
$\Gamma$ is star-shaped with just three branches. 
As in \cite{Sa}, we can take a disjoint union of subtrees $\Gamma_0$ containing the Wu set associated 
with $c$, such that the complement of $\Gamma_0$ in $\Gamma$ is a single vertex $v_0$. 
Then $\partial P(\Gamma_0 )$ is a union of the lens spaces $L_i$ and 
$P(\Gamma_0 )$ can be embedded into the interior of $P(\Gamma )$ so that $c$ extends to the spin 
structure on 
the complement $X_0 =P(\Gamma )\setminus P(\Gamma_0 )$ and on $L_i$ (we denote them by the same symbol 
$c$). Next we consider the closed $V$ manifold $\widehat X$ obtained from 
$X_0$ by attaching the cones $cL_i $ over $L_i$ and the cone $cS$ over $S$ (with orientation reversed). 
Then $c$ on $X_0$ extends naturally to the spin structure on $\widehat X$ by Proposition 5. 
Since $b_1 (\widehat X )=0$ and $b_2 (\widehat X )=1$, Theorem 1 shows that 
\[ 0=\ind D(\widehat X )=-(\sign \widehat X +\sum \delta (L_i ,c) -\delta (S,c) )/8 . \] 
Since $\sum \delta (L_i,c)=\sign P(\Gamma_0 ) - w\cdot w$ by \cite{Sa} 
and $\sign \widehat X +\sign P(\Gamma_0 ) =\sign P(\Gamma )$ by the 
additivity of the signature, 
we obtain the 
desired result. Since $w=0$ if $P(\Gamma )$ is spin, the last claim follows. 
\end{proof}

For any given spherical manifold $S$ with a spin structure $c$, we can construct a 
plumbed 4--manifold bounded by $(S,c)$ from the data of the Seifert invariants of $S$ and obtain 
the Wu set explicitly. For example, from the Seifert invariants 
$\{ (a_1 ,b_1 ), (a_2 ,b_2 ), (a_3 ,b_3 ) \}$ of $S$ 
and the data $c(g_i)$, $c(h)$ given in the list in Proposition 6, 
we can obtain anther representation of $S$ of the form 
\[ \{ (1,a), (a_1 ,b_1 '), (a_2 ,b_2 '), (a_3 ,b_3 ') \} \] 
such that $a$ is even, $a_i$ and $b_i '$ have opposite parity, and $c$ satisfies $c(g_i) =c(h)=0$
as the element of $\Hom (H_1 S^3 \setminus L ,\bz ),\bz_2 )$, where $L$ is a framed link in Figure 
\ref{Seifert} (obtained by replacing the framings $a_i /b_i$ and $0$ by $a_i /b_i '$ and $-a$ respectively). 
Then by using the continued fraction expansions of $a_i /b_i '$ by nonzero 
even numbers and by Proposition \ref{expand}, we obtain a spin plumbed 4--manifold bounded by $(S,c)$.  
This provides us an alternative method of computing $\delta (S, c)$. The details are 
omitted. 

Combining the list in Proposition \ref{delta} with the 10/8 theorem we can derive certain 
information on the intersection form of a spin 4--manifold bounded by 
a spherical 3--manifold. 

\begin{thm} \label{spherical}
Let $(S, c)$ be a spherical 3--manifold with a spin structure $c$. 
\begin{enumerate}
\item If $\delta (S,c) \ne 0$, then a connected sum of any copies of $(S,c)$ does not bound a 
$\bq$ acyclic spin 4--manifold. In particular, any $\bz_2$ homology 3--sphere $S$ with 
$\delta (S,c) \ne 0$ for a unique $c$ has infinite order in $\Theta^3_{\bz_2}$. 
\item If $|\delta (S,c)| \le 18$ and $(S,c)$ bounds a spin definite 
4--manifold $Y$, then we must have 
$\sign (Y) =\delta (S,c)$. 
\end{enumerate}
\end{thm}

\begin{proof}
The claim (1) is deduced from Proposition \ref{10/8,3}. To prove (2), we note that 
if $|\delta | <10$ then the region of $(b^- (Y), b^+ (Y))$ 
given by the two inequalities in Proposition \ref{10/8,3} does not 
contain the part with $b^+ (Y)=0$ nor $b^- (Y)=0$. If $10 \le |\delta | \le 18$, then 
the intersection of the region defined by the same inequalities and the line $b^+ (Y)=0$ or 
$b^- (Y) =0$ does not contain the point satisfying $b^+ (Y) -b^- (Y) \equiv \delta \pmod {16}$, which 
violates the condition $\sign Y\equiv \delta (S,c) \pmod {16}$. Hence we have $\sign (Y)=\delta (S,c)$. 
\end{proof}

We do not know whether a given $(S,c)$ bounds 
a definite spin 4--manifold in general, 
but in certain cases we can give such examples explicitly (see the Addendum below). 
To describe them we need some notation and results. 

\begin{figure}[ht!]
\begin{center}
\includegraphics{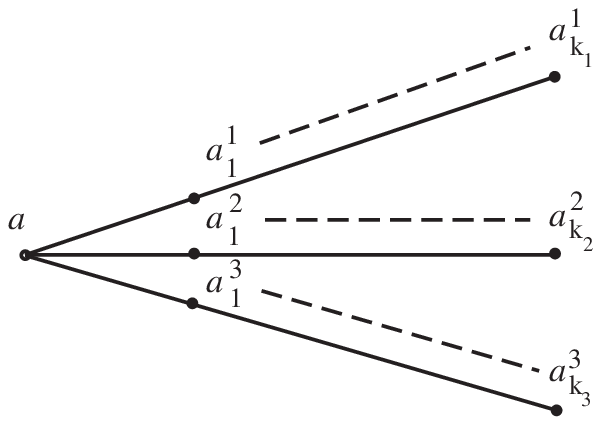}
\caption{\label{plumb} }
\end{center}
\end{figure}

\begin{nota} 
We denote the plumbed 4--manifold associated with the star-shaped diagram with 
three branches such that the weight of the central vertex is $a$ and 
the weights of the vertices of the $i$th branch are given by $(a_1^i, \dots ,a_{k_i}^i )$ 
as in Figure \ref{plumb} by 
\[ (a; a_1^1, \dots, a_{k_1}^1; a_1^2 ,\dots ,a_{k_2}^2 ; a_1^3 ,\dots ,a_{k_3}^3 ) . \] 
\end{nota}

\begin{figure}[ht!]
\begin{center}
\includegraphics{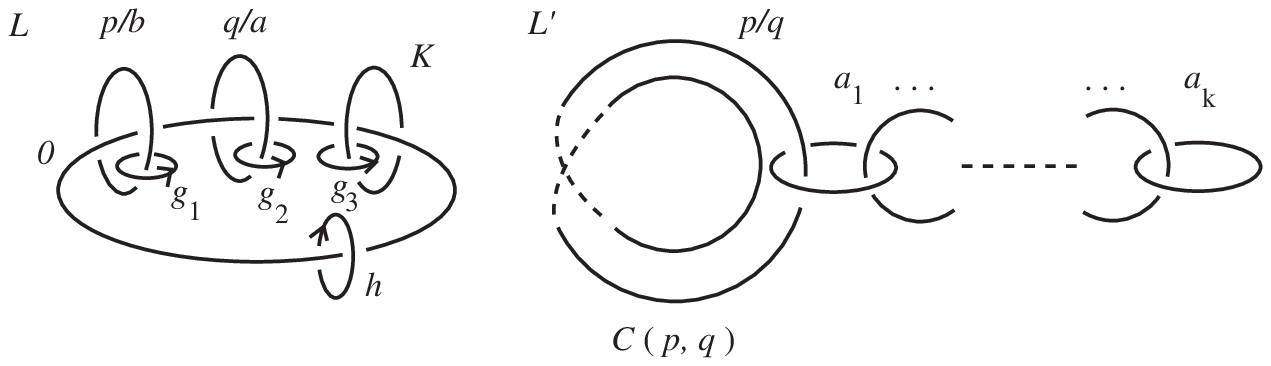}
\caption{\label{torus} }
\end{center}
\end{figure}

\begin{prop}{\rm[FS]}\qua \label{torus surgery} 
Consider a 3--manifold $M$ represented by $s/t$ surgery along a knot $K$ in a framed link $L$ in Figure 
\ref{torus}. 
Here $p$, $q$, $a$, $b$ are intergers satisfying $pa+qb=1$, and $s$ and $t$ are coprime integers 
with opposite parity. Suppose also that $M$ has a spin structure represented by 
$c\in \Hom (H_1 (S^3 \setminus L ,\bz ),\bz_2 )$ with 
$c(g_i )=w(h) =0$. Then for a 
continued fraction expansion $-t/s =[[ a_1 ,\dots ,a_k ]]$ with $a_i$ nonzero and even, 
$(M,s)$ bounds a spin 4--manifold represented by a framed link $L'$ in Figure \ref{torus}. 
Here the component of $L'$ on the left hand side is a $(p,q)$ torus knot $C(p,q)$. We denote 
$L'$ by $C(p,q) (pq; a_1 ,\dots ,a_k )$. 
\end{prop}

\begin{proof}
The knot $K$ in Figure \ref{torus} represents $C(p,q)$ in $S^3$, and 
the meridian and the preferred longitude of $K$ is given by 
$g_3$ and $h+pq g_3$ respectively. Thus $M$ 
is realized as a $pq+s/t$ surgery along $K$, and hence by framed link calculus 
we see that $M$ is also represented by $L'$, where the spin structure $c$ corresponds to 
the zero element in $\Hom (H_1 (S^3 \setminus L' ,\bz ),\bz_2 )$ by Proposition 
\ref{expand}. This proves the claim. 
\end{proof}

\begin{figure}[ht!]
\begin{center}
\includegraphics{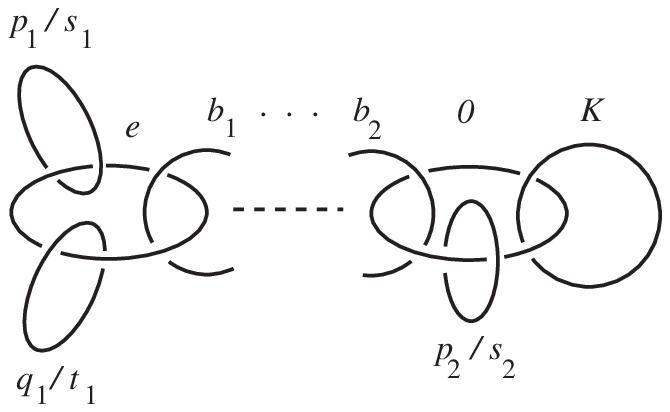}
\caption{\label{cable} }
\end{center}
\end{figure}

\begin{prop}{\rm\cite{FS}}\qua \label{iterated} 
Consider a knot $K$ in Figure \ref{cable}, where $p_1$, $q_1$, $p_2$, $s_2$, 
$s_1$, $t_1$, $e$, $b_j$ 
are integers such that 
$p_1 t_1 +q_1 s_1 =1$, $[[b_1 ,\dots ,b_s ,0, p_2/s_2 ]] =1/e$, and that 
if we put $q_2 /t_2 =[[b_s,\dots ,b_1 ]]$ for $q_2$, $t_2$ coprime, we have $p_2 t_2 +q_2 s_2 =1$. 
Then $K$ represents a knot $C(q_2 +p_2p_1 q_1 ,p_2, C(p_1 ,q_1))$ in $S^3$.  
Here we denote by $C(q,p, K)$ the cable of the knot $K$ with linking number $q$ and winding number 
$p$. 
Moreover $1/u$ surgery along $K$ in Figure 
\ref{cable} yields a $p_2 (q_2 +p_2 p_1 q_1 ) +1/u$ surgery along this 
cable knot in $S^3$, and the resulting manifold is a Seifert manifold of the form 
$\{ (1,-e), (p_1 ,s_1 ), (q_1 ,t_1 ), (r_1 ,u_1 ) \}$, where 
\[ r_1 /u_1 =[[ b_1 ,\dots ,b_s, -u, p_2/q_2 ]] .\] 
\end{prop} 

See \cite{FS} for the proof of Proposition \ref{iterated}. 
We denote by $C(u,v,C(p,q)) (a; b)$ 
the 4--manifold represented 
by a framed link with two components, which consists of $C(u,v,C(p,q))$ with framing $a$ and its 
meridian with framing $b$. 

\begin{add}
For a spherical 3--manifold $S$ with a spin structure $c$, there exists a definite spin 4--manifold 
$Y$ with $\sign Y=\delta (S,c)$ bounded by $(S,c)$ if at least $(S,c)$ satisfies one of the 
following conditions. 
\begin{enumerate}
\item 
If the negative definite 4--manifold $Y$ obtained by a minimal resolution 
of $cS$ is spin and the induced spin struction on $S$ is $c$, then $Y$ satisfies the 
above condition. 
For the cases 3--5 in Proposition \ref{delta}, the list of such $(S,c)$ and $Y$ 
is given by Table \ref{resolution}. 
\item For the cases in 3--5 in Proposition \ref{delta} other than Table 
\ref{resolution}, the minimal resolution is 
non-spin, but (S,c) in Table \ref{definite} bounds another definite spin 4--manifold $Y$ with 
$\sign Y=\delta (S,c)$. 
\end{enumerate}
\begin{table}
\caption{\label{resolution} }
\begin{tabular}{lll}
$(S,c)$ & $Y$ & $\delta$ \\
\hline 
(3-3) & $(-2k-2; -2; -2,-2; -2,-2 )$ & $-6$ \\
(4-5) & $(-2k-2; -2; -4; -2,-2 )$ & $-5$ \\
(4-13) & $(-2k-2; -2; -2,-2; -2,-2,-2 )$ & $-7$ \\
(5-5) & $(-2k-2; -2, -2,-2; -2,-2,-2,-2)$ & $-8$ 
\end{tabular}
\end{table}

\begin{table}
\caption{\label{definite} } 
\begin{tabular}{llr}
$(S,c)$ & $Y$ & $\delta$ \\
\hline
(3-5) & $C(2,-3) (-6; -2k-2, -2,-2)$ & $-4$ \\
(4-2) & $C(2,-3) (-6; -2k-2, -2, -2, -2)$ & $-5$ \\
(4-8) & $C(3,-4)(-12; 2k, -2 )$ & $-3$ \\
(4-10) & $C(2,-3) (-6; -2k-2,-4)$ & $-3$ \\
(5-1-1) & $C(2,-5) (-10; -2k-2, -2, -2)$ & $-4$ \\
(5-7) & $C(2,-3) (-6; -2k-2, -2,-2,-2,-2)$ & $-6$ \\
(3-6) with $k=-1$ & $C(2,3) (8; 2)$ & $2$ \\
(3-6) with $k=-2$ & $C(13,2, C(2,3))(26; 2)$ & $2$ \\
(4-4) with $k=-1$ & $C(2,3) (8; 2,2)$ & $3$ \\
(4-11) with $k=-1$ & $C(3,-4) (-10)$ & $-1$ \\
(4-12) with $k=-1$ & $C(2,3) (10)$ & $1$ \\
(4-14) with $k=0$ & $C(2,-3) (-2)$ & $-1$ \\
(5-2-1) with $k=-1$ & $C(2,5) (12;2)$ & 2 \\
(5-2-1) with $k=-2$ & $C(21,2, C(2,5)) (42,2)$ & $2$ \\
(5-8) with $k=-1$ & $C(2,3) (8;2,2,2)$ & $4$ \\
(5-11-(-1)) with $k=0$ & $C(2,-3) (-4;  -2)$ & $-2$ \\
\end{tabular}
\end{table}
\end{add}

\begin{figure}[ht!]
\begin{center}
\includegraphics{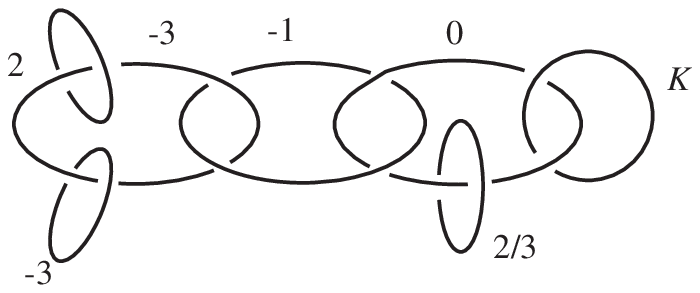}
\caption{\label{cable2} }
\end{center}
\end{figure}

\begin{proof} 
The first claim follows from Proposition \ref{wu}. The construction of $Y$ in Table 
\ref{resolution} is given 
according to the procedure explained in the paragraph after Proposition \ref{wu}. The construction of 
$Y$ for the case when the associated framed link contains a torus knot component is given 
by the procedure in Proposition \ref{torus surgery}. 
To construct $Y$ for the case (3-6) with $k=-2$, 
consider the knot $K$ in Figure \ref{cable2}. Then the $-1/2$ surgery along $K$ gives the 
Seifert manifold of type $\{ (1,3), (2,1), (3,-1), (3,-1) \}$, which is $S$ in case (3-6) with 
$k=-2$. We also note that 
$K$ in Figure \ref{cable2} gives the knot $C(1+2\cdot 6 ,2 , C(2,3))$ and $-1/2$ surgery on $K$ 
yields the $2(1+2\cdot 6) -1/2$ surgery along the cable knot. Thus 
according to Proposition \ref{iterated} the resulting 
3--manifold bounds a 4--manifold of type $C(13,2, C(2,3)) (42;2)$, 
which is spin with signature 2. It follows that this 4--manifold induces 
the (unique) spin structure on $S$ described in (3-6). 
The other cases are proved similarly and we omit the details. 
\end{proof} 

\begin{rem} \label{misc} 
\begin{enumerate}
\item If $S$ is a Poincare homology sphere (which is the case $(5-5)$ with $k=0$), the above 
result together with the classification of the unimodular forms implies that 
the intersection form of a spin negative definite 4--manifold bounded by $S$ must be $E_8$, 
which is a part of Fr\o yshov theorem [Fr] (which was also observed by Furuta). 
\item The non-unimodular definite quadratic forms over $\bz$ 
with given rank and determinant are far from unique even if the ranks are small, 
in contrast to the unimodular ones. Y. Yamada \cite{Y} pointed out that 
the Seifert manifold $\{ (2,-1), (3,-1), (5,8) \}$ ((5-2-1) with $k=-1$ in Table \ref{definite}) 
bounds three 1-connected positive definite spin 4--manifolds with intersection matrices 
$\begin{pmatrix} 2 & 1 \\ 1 & 12 \end{pmatrix}$, 
$\begin{pmatrix} 4 & 3 \\ 3 & 8 \end{pmatrix}$, 
$\begin{pmatrix} 4 & 5 \\ 5 & 12 \end{pmatrix}$ respectively, which are not congruent. 
\item $(S,c)$ with $\delta (S,c) >0$ in Table \ref{definite} bounds a negative definite (non-spin) 
4--manifold $Z$ (coming from the minimal resolution) and a positive definite spin 4--manifold $Y$, 
both of which are 1-connected. Then 
$Z\cup (-Y)$ is a closed 1-connected negative definite 4--manifold, which is homeomorphic to 
a connected sum of $\overline {\mathbf{CP}}^2$'s by Donaldson and Freedman's theorem. 
Y. Yamada [Y] observed that such a manifold appeared in Table \ref{definite} other than 
the cases (3-6) and (5-2-1) (with $k=-2$) is diffeomorphic 
to $\sharp \overline{\mathbf{CP}}^2$. 
\end{enumerate}
\end{rem}

\begin{rem}\label{BL}
In \cite{BL} Bohr and Lee considered two invariants $m(\Sigma )$ and $\overline m (\Sigma )$ 
for a $\bz_2$ homology 3--sphere $\Sigma$. 
Here $m(\Sigma )$ is defined as the maximum of 
$5\sign X/4 -b_2 (X) $, while $\overline m (\Sigma )$ is defined as 
the minimum of $5\sign X/4 +b_2 (X)$, where 
$X$ ranges over all spin 4--manifolds with $\partial X =\Sigma$. 
They proved that if $m(\Sigma )>0$ or $\overline m (\Sigma )<0$ 
(or $m(\Sigma )=0$ and $R(\Sigma )\ne 0$) then $\Sigma$ has infinite order in $\Theta^3_{\bz_2 }$ by 
using the $10/8$ theorem for spin 4--manifolds. 
For example, consider $\Sigma =\{ (2,-1), (3,-1), (3, -6k+2 ) \}$ with $k<0$ 
((3-4) in Proposition \ref{delta}). Then 
as is remarked after Proposition 8, we can see that $\Sigma$ bounds 
a spin plumbed 4--manifold $P(\Gamma )$ with $b^+ (P(\Gamma ) )=5$ and $b^- (P(\Gamma )) =1$. It follows that 
$-1 \le m (\Sigma )  \le \overline m (\Sigma ) \le 11$. 
If $\Sigma$ would bound a spin definite 4--manifold $Y$ (we do not know whether this is the case), then 
we could see that $m(\Sigma )>0$ and recover the claim (1) in Theorem 2 for this case. 
M. Furuta \cite{Fu2} pointed out that if we extend 
the above definitions of $m(\Sigma )$ and $\overline m (\Sigma )$ 
by replacing $5\sign X/ 4 \pm b_2 (X)$ by $5\sign X /4 +\delta (X)/4 \pm b_2 (X)$, where 
$X$ ranges over all spin $V$ 4--manifolds with $\partial X =\Sigma$, 
then the same conclusion in \cite{BL} is obtained by virtue of the $V$ version of the 
$10/8$ theorem. 
\end{rem}

Finally we give an application of Theorem \ref{spherical} 
to embeddings of $\mathbf {RP}^2$ into 4--manifolds. 
The following theorem generalizes the result in \cite{L} in the case when the embedded 
$\mathbf{RP}^2$ is a characteristic surface. 

\begin{thm} \label{projective}
Let $X$ be a closed smooth 4--manifold $X$ with $H_1 (X, \bz ) =0$. Suppose that 
there exists 
a smoothly embedded real projective plane $F$ in $X$ with $PD [F] \pmod 2 =w_2 (X)$. 
Denote by $e(\nu )$ the normal Euler number of 
the normal bundle of the embedding $F \subset X$. Then 
$\sign X -e(\nu ) \equiv \pm 2 \pmod {16}$. Furthermore 
\begin{enumerate}
\item If $\sign X -e(\nu )\equiv 2\epsilon \pmod {16}$ with 
$\epsilon =\pm 1$, then either $e(\nu )+2\epsilon =\sign X$ or 
\[ 8(1-b^- (X)) +\sign X \le e(\nu ) +2\epsilon \le 8(b^+ (X) -1)+\sign X .\] 
\item If both $b^- (X) <3$ and $b^+ (X) <3$, then $e(\nu ) =\sign X -2\epsilon$. 
\end{enumerate}
\end{thm}

\begin{proof} 
Let $[F]$ be the element of $H_2 (X, \bz_2 )$ represented by $F$. 
First suppose that there exists an element $y\in H_2 (X, \bz_2 )$ with $[F] \cdot y \equiv 1 \pmod 2$, and 
$e(\nu ) \ne 0$. Put $n=e(\nu )$. 
Let $X_0$ be the complement of the tubular neighborhood $N(F)$ of $F$ in $X$. 
Then $\partial N(F) =-\partial X_0$ is the twisted $S^1$ bundle over $\mathbf {RP}^2$ with normal 
Euler number $n$, which is diffeomorphic to a Seifert manifold over $S^2 (2,2,|n| )$ (which we denote by 
$S$) with Seifert invariants $\{ (2,1), (2,-1), (|n|,\sgn n) \}$. 
We fix the correspondence between $\partial N(F)$ and $N$ as follows. 
Denote by $g_i$, $h$ be the curves in the framed link picture $L$ of $S$ 
associated with the above Seifert 
invariants as in Figure \ref{Seifert}. 
Also denote by $Q$, $H$ be (one of) the cross section of the curve generating 
$H_1 (\mathbf {RP}^2 ,\bz_2 )$ and the fiber of the $S^1$ bundle $\partial N(F)$. 
Then we have a diffeomorphism between 
$S$ and $\partial N(F)$ so that 
\[ Q=g_2 ,\quad H=g_1 +g_2 =-g_3, \quad 2Q =h \] 
in the first homology group. Considering the 
exact sequence of the homology groups for the pair $(X, X_0 )$, we see that $H_2 (X_0 ,\bz_2 )$ 
is the set of $x\in H_2 (X, \bz_2 )$ with $x\cdot [F]  \equiv \langle w_2 (X) ,x \rangle
 \equiv 0 \pmod 2$. Hence 
$X_0$ admits a spin structure $c$ (which is unique 
since $H_1 (X_0 ,\bz _2 ) =0$ by the above assumption). 
Since the spin structure induced on $S=-\partial X_0$ from $c$ extends 
uniquely to that on the cone $cS$ over $S$ (Proposition \ref{extend}), 
$c$ extends uniquely to the spin structure on 
$\widehat X =cS \cup X_0$ (which we also denote by $c$). 
Since $H_3 (X, X_0 ,\bz )=H^1 (F, \bz ) =0$ and 
$H_2 (X, X_0 ,\bz ) =H^2 (F, \bz ) =\bz_2$, we have $H_2 (X_0 ,\bq ) =H_2 (X, \bq)$ and hence 
$b^{\pm} (\widehat X)=b^{\pm} (X)$ and $\sign (\widehat X) =\sign (X)$.  
Note that since 
$\langle w_2 (X) ,y \rangle \equiv [F] \cdot y \equiv 1 \pmod 2$, 
the spin structure restricted on $H$ does not 
extend to that on the disk fiber of $N(F)$. Under the above correspondence, this implies that if $c$ 
restricted on $S$ is represented by a 
homomorphism from $H_1 (S^3 \setminus L ,\bz )$ to $\bz_2$ 
for a framed link $L$ as in Definition \ref{spin}, then $c(H) 
\equiv c(q_3 ) \equiv c(q_1 )+c (q_2 ) \equiv 1 \pmod 2$. 

\subsection*{The case when $n>1$} 
Consider the representation of $S$ by $\{ (2,1), (2,1), (n,1-n ) \}$. 
If we denote the curves associated with 
the corresponding framed link picture by $g_i '$ and $h'$, then the correspondence between them and the 
original curves is given by 
\[ g_1 '=g_1 ,\quad g_2 ' =g_2 -h, \quad g_3 ' =g_3 +h, \quad h' =h .\] 
Now we check $\delta (S,c)$ according to the list in Proposition 
\ref{delta} (note that we have always $c(h)=c(h') =0$). 
\subsubsection*{The case when $(c(g_1 '), c(g_2 '), c(g_3 ' ) ) =(0,0,0)$} 
Under the above correspondence we have $(c(g_1 ), c(g_2 ) , c(g_3 ))=(0,1,1)$ 
and hence $c(Q)=c(H)=1$. Then since $n/(n-1) =\overbrace{[[2,\dots,2]]}^{n-1}$, (2-1) 
(if $n$ is odd) or (2-6) (if $n$ is even) in Proposition \ref{delta} 
shows that 
\[ \delta (S,c)=\sigma (n,1-n,-1)=-\sigma (n,n-1,-1)=\sigma (n-1 ,n,-1) +1 =-n+2. \]  
It follows that $\ind D(\widehat X )=-(\sign (\widehat X )+\delta (S,c))/8 =-(\sign X -n+2)/8$. 
Thus Proposition \ref{10/8,3} shows that $n \equiv \sign X+2 \pmod {16}$, and either 
$n =\sign X +2$ or 
\[ 8(1-b^- (X)) +\sign X  \le n-2 \le 8(b^+ (X) -1) +\sign X .\] 
\subsubsection*{The case when $(c(g_1 ' ), c(g_2 ' ) ,c(g_3 '))=(1,1,0)$} In this case 
$(c(g_1 ), c(g_2 ), c(g_3 )) =(1,0,1)$ and hence $c(Q)=0$, $c(H)=1$. Then (2-2) or (2-7) 
in Proposition \ref{delta} implies that 
\[ \delta (S,c)=\sigma (n,1-n,-1) -4 =-n-2 .\] 
Thus $n \equiv \sign X -2 \pmod {16}$, and either $n=\sign X -2$ or 
\[ 8(1-b^- (X)) +\sign X \le n+2 \le 8(b^+ (X) -1) +\sign X .\] 
\subsubsection*{The case when $(c (g_1 '),c(g_2 ') ,c(g_3 '))=(\epsilon ,1-\epsilon ,1)$ and $n$ is even} 
In this case we have $(c (g_1 ) ,c(g_2 ), c(g_3 ))=(\epsilon ,\epsilon, 0)$ and hence $c(H)=0$. But this violates the 
above condition, and hence this case cannot occur. 

\subsection*{The case when $n<-1$} 
Reversing the orientation of $X$ we have an embedding $F \subset -X$ whose normal Euler number is $-n$. 
If we consider $(-\widehat X ,-c)$ in place of $(\widehat X ,c)$, 
we have $b^+ (-\widehat X )=b^- (\widehat X) =b^- (X)$, 
$b^- (-\widehat X )=b^+ (\widehat X )=b^+ (X)$, $\sign (-\widehat X )=-\sign \widehat X =-\sign X$, and 
$\delta (-S, -c)=-\delta (S,c)$. Hence we derive the same result from the case when $n>1$ by applying 
Proposition \ref{10/8,3} to $-\widehat X$.  

Next we consider the general case. Consider the internal connected sum of $F \subset X$ and $k$ copies of 
the standard embedding $\mathbf {CP}^1 \subset \mathbf {CP}^2$ for some $k$ to obtain another 
embedding $\widetilde F \subset \widetilde X$ of $\mathbf {RP}^2$, 
where $\widetilde F =F\sharp k \mathbf{CP}^1$ 
and $\widetilde X=X\sharp k\mathbf{CP}^2$. Then $ PD [\widetilde F ] \pmod 2 =w_2 (\widetilde X )$ and 
there exists an element $y \in H_2 (\widetilde X ,\bz_2 )$ with $y\cdot {\widetilde F } \equiv 1 \pmod 2$ (for 
example, choose a copy of $\mathbf {CP}^1$ in one $\mathbf {CP}^2$ summand as $y$). Moreover 
the normal Euler number $e(\widetilde {\nu} )$ of the embedding $\widetilde F \subset \widetilde X$ is 
$e(\nu ) +k$ (which is greater than one for some $k$), 
$\sign (\widetilde X )=\sign X +k$, and $b^- (\widetilde X ) =
b^- (X)$. Thus applying the above result to $\widetilde F \subset \widetilde X$, 
we have 
$\sign X -e(\nu ) =\sign (\widetilde X ) -e(\widetilde {\nu} ) \equiv \pm 2 \pmod {16}$, and obtain the 
inequality on the left hand side in (1). If we consider the embedding $\widetilde F \subset 
\widetilde X$ obtained by 
the internal connected sum of 
$F \subset X$ and $k$ copies of the standard embedding 
$\overline {\mathbf {CP}}^1 \subset \overline {\mathbf {CP}}^2$, 
we have $\widetilde F =F \sharp k \overline{\mathbf {CP}}^1$, 
$\widetilde X =X\sharp k\overline{\mathbf {CP}}^2$, 
the normal Euler number of $\widetilde F \subset \widetilde X$ is $e(\nu ) -k$ 
(which is less than $-1$ for some $k$), 
$\sign \widetilde X =\sign X -k$, and $b^+ (\widetilde X )=b^+ (X)$. 
We also have $y \in H_2 (\widetilde X ,\bz_2 )$ with $y \cdot \widetilde F \equiv 1 \pmod 2$. 
Thus applying the above result to this embedding 
we obtain the inequality on the right hand side in (1). 
To see (2) suppose that $b^+ (X) <3$ and $b^- (X) <3$. Then $\sign X -16< 8(1-b^- (X)) +\sign X$ and 
$8(b^+ (X)-1) +\sign X <\sign X +16$. Since $e(\nu )+2\epsilon \equiv \sign X \pmod {16}$, 
the above inequalities do not hold unless 
$e(\nu ) + 2\epsilon =\sign X$. This proves (2). 
\end{proof}

\begin{rem}
The claim $e(\nu ) -\sign X \equiv \pm 2 \pmod {16}$ is also deduced from Guillou and Marin's theorem 
\cite{GM}, \cite{M}. 
We note that for some $X$ both of the cases when $\epsilon =\pm 1$ in 
(1) occur. For example, when $X=k\mathbf {CP}^2$, consider the 
connected sum of $\mathbf {RP}^2 \subset S^4$ with normal Euler class $\pm 2$ and the $k$ copies of 
$\mathbf {CP}^1 \subset \mathbf{CP}^2$ (\cite{L}). 
\end{rem}

\begin{rem}
Acosta \cite{A} obtained the estimate on the self-intersection number of a characteristic element $x$ of $X$ 
that is realized by a smoothly embedded 2--sphere by considering the $10/8$ theorem for $V$ 4--manifolds with 
$cL(p, \pm 1  )$ type singularities. In our terminology, the result is derived as follows. 
Suppose that 
$X$ is non-spin and $x$ is realized by an embedded 2--sphere $F$ with $x\cdot x =n >0$. Then the 
complement $X_0$ of the tubular neighborhood $N(F)$ in $X$ has a spin structure $c$, 
which does not extend to that 
on the disk fiber of $N(F)$. 
Since $N(F)$ is represented by the $n$ surgery along the trivial knot $O$ in $S^3$, 
this implies that the induced spin structure on $\partial N(F) =L(n,-1)$ corresponds to the homomorphism 
from $H_1 (S^3 \setminus O ,\bz )$ (generated by the meridian $\mu$ of $O$) to $\bz_2 $ with $c(\mu ) =1$. 
Consider the $V$ manifold $\widehat X =c L(n,-1) \cup X_0$ with spin structure $c$, 
which is an extension of the 
original one. Then 
according to \cite{FFU} Proposition 3, the contribution $\delta (L(n ,-1), c)$ 
to $\ind D(\widehat X )$ equals 
$\sigma (-1, n, (-1)^{c(\mu )-1} )=\sigma (-1, n, 1) =\sigma (n-1, n, -1) =-(n-1)$. 
Thus applying the $10/8$ theorem to $\widehat X$ we obtain the 
inequality in \cite{A}. The general case is proved by a similar argument 
as in the proof of Proposition \ref{projective} (\cite{A}). 
\end{rem}

\Addresses\recd

\end{document}